\documentclass[12pt]{article}
\usepackage{amsmath,amsthm,amsfonts,amscd,amssymb} 
\usepackage{verbatim}   
\usepackage{eucal} 
\usepackage{longtable} 
\usepackage{makeidx}

\usepackage[T1]{fontenc}
\usepackage[utf8]{inputenc}
\usepackage{authblk}

\newcommand{\field}[1]{\mathbb{#1}}
\newcommand{\califas}[1]{\mathcal{#1}}

\newcommand{\R}{\field{R}}
\newcommand{\C}{\field{C}}

\newcommand{\Hc}{\califas{H}}

\newfont{\cyrr}{wncyr10}

\newcommand{\be}{\begin{equation}}
\newcommand{\bel}[1]{\begin{equation}\label{#1}}
\newcommand{\ee}{\end{equation}}

\newcommand{\ba}{\begin{eqnarray}}
\newcommand{\ea}{\end{eqnarray}}

\newcommand{\bas}{\begin{eqnarray*}}
\newcommand{\eas}{\end{eqnarray*}}

\newcommand{\thm}{\begin{theorem}}
\newcommand{\thml}[1]{\begin{theorem}\label{#1}}
\newcommand{\mht}{\end{theorem}}
\newcommand{\cnj}{\begin{conjecture}}
\newcommand{\cnjl}[1]{\begin{conjecture}\label{#1}}
\newcommand{\jnc}{\end{conjecture}}
\newcommand{\dfn}{\begin{definition}}
\newcommand{\dfnl}[1]{\begin{definition}\label{#1}}
\newcommand{\nfd}{\end{definition}}
\newcommand{\pro}{\begin{proposition}}
\newcommand{\prol}[1]{\begin{proposition}\label{#1}}
\newcommand{\orp}{\end{proposition}}
\newcommand{\crl}{\begin{corollary}}
\newcommand{\crll}[1]{\begin{corollary}\label{#1}}
\newcommand{\lrc}{\end{corollary}}
\newcommand{\lmm}{\begin{lemma}}
\newcommand{\lmml}[1]{\begin{lemma}\label{#1}}
\newcommand{\mml}{\end{lemma}}
\newcommand{\prf}{\begin{proof}}
\newcommand{\prfl}[1]{\begin{proof}\label{#1}}
\newcommand{\frp}{\end{proof}}
\newcommand{\axi}{\begin{axiom}}
\newcommand{\axil}[1]{\begin{axiom}\label{#1}}
\newcommand{\ixa}{\end{axiom}}
\newcommand{\rmk}{\begin{remark}}
\newcommand{\rmkl}[1]{\begin{remark}\label{#1}}
\newcommand{\kmr}{\end{remark}}
\newcommand{\exa}{\begin{example}}
\newcommand{\exal}[1]{\begin{example}\label{#1}}
\newcommand{\axe}{\end{example}}
\newcommand{\alg}{\begin{algorithm}}
\newcommand{\gla}{\end{algorithm}}
\newcommand{\nte}{\begin{note}}
\newcommand{\ntel}[1]{\begin{note}\label{#1}}
\newcommand{\etn}{\end{note}}
\newcommand{\app}{\begin{application}}
\newcommand{\appl}[1]{\begin{application}\label{#1}}
\newcommand{\ppa}{\end{application}}
\newcommand{\mat}{\begin{matrix}}
\newcommand{\tam}{\end{matrix}}
\newcommand{\smm}{\begin{summary}}
\newcommand{\mms}{\end{summary}}

\newcommand{\teo}{\begin{teorema}}
\newcommand{\teol}[1]{\begin{teorema}\label{#1}}
\newcommand{\oet}{\end{teorema}}
\newcommand{\df}{\begin{definicion}}
\newcommand{\dfl}[1]{\begin{definicion}\label{#1}}
\newcommand{\fd}{\end{definicion}}
\newcommand{\por}{\begin{proposicion}}
\newcommand{\porl}[1]{\begin{proposicion}\label{#1}}
\newcommand{\rop}{\end{proposicion}}
\newcommand{\cor}{\begin{corolario}}
\newcommand{\corl}[1]{\begin{corolario}\label{#1}}
\newcommand{\roc}{\end{corolario}}
\newcommand{\lem}{\begin{lema}}
\newcommand{\leml}[1]{\begin{lema}\label{#1}}
\newcommand{\mel}{\end{lema}}
\newcommand{\pru}{\begin{prueba}}
\newcommand{\prul}[1]{\begin{prueba}\label{#1}}
\newcommand{\urp}{\end{prueba}}
\newcommand{\axa}{\begin{axioma}}
\newcommand{\axal}[1]{\begin{axioma}\label{#1}}
\newcommand{\xa}{\end{axioma}}
\newcommand{\nta}{\begin{nota}}
\newcommand{\ntal}[1]{\begin{nota}\label{#1}}
\newcommand{\atn}{\end{nota}}
\newcommand{\eje}{\begin{ejemplo}}
\newcommand{\ejel}[1]{\begin{ejemplo}\label{#1}}
\newcommand{\je}{\end{ejemplo}}
\newcommand{\api}{\begin{aplicacion}}
\newcommand{\apil}[1]{\begin{aplicacion}\label{#1}}
\newcommand{\ipa}{\end{aplicacion}}

\newcommand{\ali}{\begin{align}}
\newcommand{\ila}{\end{align}}
\newcommand{\enu}{\begin{enumerate}}
\newcommand{\une}{\end{enumerate}}
\newcommand{\arr}{\begin{array}}
\newcommand{\rra}{\end{array}}
\newcommand{\eqa}{\begin{eqnarray}}
\newcommand{\aqe}{\end{eqnarray}}
\newcommand{\equ}{\begin{equation}}
\newcommand{\uqe}{\end{equation}}
\newcommand{\subq}{\begin{subequations}}
\newcommand{\qbus}{\end{subequations}}

\newtheorem{theorem}{Theorem}[section]
\newtheorem{lemma}[theorem]{Lemma}
\newtheorem{proposition}[theorem]{Proposition}
\newtheorem{corollary}[theorem]{Corollary}
\newtheorem{conjecture}[theorem]{Conjecture}

\theoremstyle{definition}
\newtheorem{definition}{Definition}[section]
\newtheorem{axiom}[definition]{Axiom}

\theoremstyle{remark}
\newtheorem{remark}{Remark}[section]
\newtheorem{example}[remark]{Example}
\newtheorem{application}[remark]{Application}
\newtheorem{algorithm}[remark]{Algorithm}
\newtheorem{note}[remark]{Note}
\newtheorem{summary}[remark]{Summary}

\newtheorem{teorema}[theorem]{Teorema}
\newtheorem{lema}[theorem]{Lema}
\newtheorem{proposicion}[theorem]{Proposici\'on}
\newtheorem{corolario}[theorem]{Corolario}

\newtheorem{definicion}[definition]{Definici\'on}
\newtheorem{axioma}[definition]{Axioma}

\newtheorem{nota}[remark]{Nota}
\newtheorem{ejemplo}[remark]{Ejemplo}
\newtheorem{aplicacion}[remark]{Aplicaci\'on}

\begin{document}

\title{Integrable Systems and Geometry of Riemann Surfaces}

\author[$\dagger$]{ J.A. ESP\'INOLA-ROCHA}

\affil[$\dagger$]{{\small \it Departamento de Ciencias B\'asicas. 

Universidad Aut\'onoma Metropolitana-Azcapotzalco.

Azcapotzalco, 
CDMX,
M\'exico, }
{\small jaer.azc.uam.mx@gmail.com}}

\author[$\ddagger$]{F.X. PORTILLO-BOBADILLA}

\affil[$\ddagger$]{{\small \it Universidad Aut\'onoma de la Ciudad de M\'exico.
 Iztapalapa, 
CDMX,  M\'exico,}
{\small francisco.portillo@uacm.edu.mx}}

\maketitle

\begin{abstract}
We give a self-contained introduction to the relations between Integrable Systems and the Geometry of Riemann Surfaces. We start from a historical introduction to the topic of integrable systems. Afterwards, we study the polynomial solutions to the stationary KdV equation, giving concrete examples of the computations of the solutions on small degree ($N=0$, $1$, $2$.) We discused the geometry of the so called {\it square eigen-functions} for those cases. Later,
we present the solutions in the general case, presenting a formula for the hyper-elliptic curve of genus $N$ that parametrizes the solutions, the $N$-solitons. We relate also our equations to the Lax hierarchy, and analyze the time evolution for the KdV.
Using the scalar operator for the NLS equation, computed by Kamchatnov, Krankel and Umarov \cite{Kamchatnov_Kraenkel_Umarov_1, Kamchatnov_Kraenkel, Kamchatnov_Kraenkel_Umarov_2}, we analyse such an equation following a similar approach as the one we used for the KdV equation in the
first part of the manuscript. We present some original results obtained in that case.

 \end{abstract}


\newpage

\tableofcontents

\newpage

\section*{Introduction}

Many equations of the mathematical-physics are obiquous and appear in many applications and in the physics of many nature phenomena. 
To understand nature, more sophisticated models are needed. However, {\it integrable
systems} are good minimal models to describe and understand more complicated
phenomena in nature, such as wave propagation. 
 This is why we are interested in  Integrable  Systems, which can be solved by means of the inverse scattering transform  (IST).
 
 Among the most notable completely integrable partial differential equations we find the Korteweg-deVries (KdV), the Non-linear Schr\"odinger (NLS), the Derivative Non-linear Schr\"odinger (dNLS) and 
 the sine-Gordon (sG) equations. 
 
In the early 1830s, the scottish engineer Scott Russell observed a {\it large solitary elevation} in a narrow Edinburgh channel, now called the {\it Scott Russell's Aqueduct}. Russell observed the phenomenon years before the equation that describes this wave motion were discovered,  and certainty he did not have any idea of the impact that his observation and experiments would have in physics, mathematics and modern technology  \cite{scott}.
 
 Namely, the Korteweg-deVries equation is given, in its standard form, by \cite{KdV_seminal}
 \[
  \frac{\partial q}{\partial t} + 6q\frac{\partial q}{\partial x} + \frac{\partial^3 q}{\partial x^3}  =   0,
  \]
which describes the profile of a wave propagating in shallow water.  The {\em one-soliton} solution derived by D.J. Korteweg and G. deVries (1895) is given by:
\[
q(x,t) = \frac{c}{2}\text{sech}^2\left[  \frac{\sqrt{c}}{2} \left(x - ct - x_0  \right) \right].
\]
This is a solution in the infinite line, with infinite period. Korteweg and deVries also found the traveling wave solution for the periodic case, and they found the solution expressed as an elliptic function, 
\begin{equation}\label{cnoidal_sol_n}
\begin{split}
q(x,t) &= f(x - ct) \\ &= -f_2 + (f_2 - f_3) \text{cn}^2\left[ (x - ct - x_0)   \left.   \sqrt{\left\{ \frac{f_1 - f_3}{2} \right\}} \;\; \right| m \right]
\end{split}
\end{equation}
where $c$ is the propagation velocity, and $f_1$, $f_2$ and $f_3$ are  roots of the polynomial  $F(f) := -f^3+\frac{c}{2}f^2+Af+B$, where $A$, $B$ are constants of integration. (See
Section \ref{part1}.\ref{sect:kdv_torus}.)

In part \ref{part1} of the manuscript, we consider the solutions of the periodic (in $x$) KdV equation as worked by Dubrovin, Novikov, Its, Krichever and Matveev \cite{dubrovin_1, krichever_1, krichever_novikov, dubrovin_matveev_novikov}, but following Flaschka and Newell \cite{flaschka_1, newell}. 
The  analysis is performed through the scalar eigenvalue problem associated to the KdV equation, namely, the stationary Schr\"odinger equation with potential $q = q(x)$, where $q(x)$ is a periodic solution of the stationary KdV equation  \cite{flaschka_1, newell}. The case for the non-stationary solution is briefly discussed for the KdV equation in section \ref{sec3}.

In part \ref{part2},  we develop the theory studied in part \ref{part1}, but applied to the NLS equation. As opposed to the standard matrix approach of the Zakharov-Shabat (ZS) operator, we write this operator in scalar form as Kamchatnov, Kraenkel and Umarov in \cite{Kamchatnov_Kraenkel_Umarov_1, Kamchatnov_Kraenkel, Kamchatnov_Kraenkel_Umarov_2}, where the authors express some matrix eigenvalue problems (associated to Lax's operators) as a scalar eigenvalue problem, and proceed as in part \ref{part1}. We recover the infinite set of NLS fluxes and constants of motion as a recursion formula as it was done in the literature \cite{zajarov} and, at each step of recursion, we find two fluxes. One of them corresponds to the constants of motion, and the other one to the Lax hierarchy \cite{espinola_portillo}.
We also write the Riemann hyperelliptic curve associated to this equation.




In part \ref{part3} of this work, we include three apendices
which are needed for the results in the parts \ref{part1} and \ref{part2}. 
The appendix \ref{apnA} is the proof of one of the main results Theorem \ref{thm2}.
The appendix \ref{apnB} is the proof of a well known identity of symetric functions,
which is usually used in the literature, but here you can find a novel and elementary proof of that identity. Appendix \ref{apnC} contains the computations of
the derivatives of $E$, which are used to simplify expressions in the NLS computations in part \ref{part2}.



\part{The KdV equation}\label{part1}

\section{KdV equation and the Complex Torus}\label{sect:kdv_torus}

Consider the KdV equation $q_t+6qq_x+q_{xxx}=0$ with the extra assumption that
the solution is a fixed wave that mantains its shape $q(x,t)=f(x-ct)$ (i.e. it is a travelling wave). Set $s=x-ct$.
Thus, the KdV equation can be expressed just in terms of $f$ and its derivatives and it
becomes $-cf'+6ff'+f'''=0$.  Integrating, we obtain $-cf+3f^2+f''=A$, where $A$ is a constant.
Now, multiplying by $f'$ and integrating again, it follows the equation $-\frac{c}{2}f^2+f^3+\frac{1}{2}(f')^2=Af+B$. 
The last equation corresponds to a Complex Torus or an Elliptic Curve $E$ with coordinates $(f,f')\in \C ^2$.


Using that $f'=\frac{df}{ds}=\sqrt{2(-f^3+\frac{c}{2}f^2+Af+B)}$ and integrating we have
that $$x-ct=s=\int_{P_0}^{P} \frac{df}{\sqrt{2(-f^3+\frac{c}{2}f^2+Af+B)}},$$ with $P$ a moving point and $P_0$ a chosen initial point in $E$, and the inverse of this function is given as in equation (\ref{cnoidal_sol_n}).

\section{The Lax pair for KdV equation}\label{sect:kdv_lax_pair}

As it is well known from the literature, or the introductory texts \cite{newell, novikov, drazin}, the Korteweg-deVries equation
\begin{equation}\label{kdv}
  \frac{\partial q}{\partial t} + 6q\frac{\partial q}{\partial x} + \frac{\partial^3 q}{\partial x^3}  =   0,
\end{equation}
is a completely integrable equation \cite{zakharov_faddeev}, and it is given by the compatibility condition of the two scalar differential equations \cite{ggkm}:  
\begin{eqnarray}
\widehat{\cal{L}}y  &=& \lambda y,      \label{e_val_problem}  \\
\frac{\partial y}{\partial t}  &=& \widehat{\cal{P}}y,     \label{t_flow} 
\end{eqnarray}
where $\widehat{\cal{L}}$ and $\widehat{\cal{P}}$ is a pair of differential operators,
\begin{eqnarray}
\widehat{\cal{L}} &=&  \frac{\partial^2}{\partial x^2} + q(x,t),      \label{L_operator}  \\
\widehat{\cal{P}}  &=& 4 \frac{\partial^3}{\partial x^3} + 6q(x,t) \frac{\partial}{\partial x} + 3q_x(x,t),    \label{P_operator} 
\end{eqnarray}
and $\lambda$ is the spectral parameter. This pair of operators is known as the {\em "Lax's pair"} for the KdV equation. They were called this way in honor to P. Lax after his discovering in 1968 \cite{lax_1}.

Equation (\ref{e_val_problem}) is an eigenvalue problem for the function $y(x,t)$, which turns to be (Miraculously!) the stationary Schr\"odinger equation of Quantum Mechanics. Equation (\ref{t_flow}) is the evolution in time equation, or the $t$-flow equation.

Now, $q(x,t)$ in equation (\ref{L_operator}) is the solution to the KdV equation (\ref{kdv}) if and only if 
\[
\frac{d\lambda }{dt} = 0, \makebox[1in]{ \cite{ggkm, lax_1}}
\]

We then say that the KdV equation is an {\em isospectral flow}.

\section{Traveling wave solutions of the KdV equation}\label{sect:traveling_soln_s_kdv}

In this section,  we consider stationary solutions of the KdV equation. We will study first the infinite  and  finite periodic solutions in traveling wave form. Then, we will consider the "{\em 0-soliton}", "{\em 1-soliton}", and the "{\em 2-soliton}" solutions of the periodic KdV. 
Afterwards, we will consider time-dependent solutions. All the computations here are following \cite{flaschka_1, newell}, from the original work of  \cite{dubrovin_1, krichever_1, krichever_novikov}.

In the original paper by Korteweg and deVries \cite{KdV_seminal} the authors construct the equation for waves in shallow water, and also found the periodic and infinite-periodic solutions we are going to show here. They proceeded as follows, which is also a standard computation \cite{drazin}. Consider the KdV equation (\ref{kdv}), and consider the change of variable $u(x,t)=-q(x,t)$ to obtain the equation:
\equ
u_t-6uu_x+u_{xxx}=0
\uqe
Consider thus, the
solutions of the form
\[
u(x,t) = f(\xi),
\]
donde $\xi = x-ct$, with $c$ is a constant and it is the velocity of propagation. We get,
\[
-c f' - 6ff' + f''' = 0. 
\]
Integrating twice,
\begin{equation} \label{phase_space_eqn}
\frac{1}{2}(f')^2  = \frac{c}{2} f^2 + f^3 + Af + B, 
\end{equation}
where $A$ and $B$ are constants of integration. Assuming $f, f', f'' \to 0$ as $|\xi| \to \infty$, we get $A=B=0$ and 
\[
\frac{1}{2}(f')^2  = \frac{c}{2} f^2 + f^3, 
\]
which is a separable equation 
\[
\frac{d\xi}{df} = \pm \frac{1}{f\sqrt{2f + c}},
\]
assuming that $2f + c > 0$ to get real solutions. Integrating, we obtain the {\em "1-soliton"} solution 
\[
q(x,t)=-u(x,t) = \frac{1}{2}c \:\:\text{sech}^2\left[  \frac{1}{2}\sqrt{c}(x - ct - x_0)  \right],
\]
which is a "bell-shaped" profile, as observed by Scott Russell in 1835. 

If $A\neq 0, B\neq 0$, we have equation (\ref{phase_space_eqn}), which is also separable, to get
\begin{equation} \label{separable_eqn}
\frac{d\xi}{df} = \pm \frac{1}{ \sqrt{2F(f)}},
\end{equation}
with $F(f)  = f^3  + \frac{c}{2} f^2 + Af + B$ and $F(f) \neq  0$, with solution
\[
\xi = \xi_0   \pm \int_{f_0}^f  \frac{1}{ \sqrt{2F(\hat f)}}d\hat f,
\]
with $f_0 = f(\xi_0)$, where $\xi_0$ is an initial position for the limit of integration. This is an elliptic integral.

Now, if we consider that $F(f)=0$ has three distict real roots, order as
$f_3< f_2< f_1$. And, if we take as our initial limit of integration $f=f_3(\xi_3)=f_3$, we obtain:
\begin{equation}\label{xi_soln_in_f}
\xi = \xi_3   \pm \int_{f_3}^f  \frac{1}{ \sqrt{2(\hat f - f_1)(\hat f - f_2)(\hat f - f_3)}}d\hat f,
\end{equation}

Under the change of variables,
\[
\hat f = f_3 + (f_2 - f_3)\sin^2\theta, \qquad (f_2 - f_3 \neq 0),
\]
we obtain
\begin{equation}\label{xi_soln}
\xi(\phi) = \xi_3   \pm  \sqrt{ \frac{2}{ f_1 - f_3} }  \int_{0}^{\phi}  \frac{1}{ \sqrt{1 - m \sin^2\theta}}d\theta,
\end{equation}
with $m = (f_2 - f_3)/(f_1-f_3)$ and
\begin{equation}\label{pre_cnoidal_soln}
f = f_2 - (f_2 - f_3)\cos^2\phi
\end{equation}
relates the upper limits of the integrals (\ref{xi_soln_in_f}) and (\ref{xi_soln}).

Now, if in the integral 
$$ \int_{0}^{\phi}  \frac{1}{ \sqrt{1 - m \sin^2\theta}}d\theta,$$
we invert
$\xi = \xi(\phi)$ to $\phi = \phi(\xi)$, and we define the {\em "cnoidal-cosine"},
\[
\text{cn}(\xi | m) = \cos(\phi(\xi)). 
\]
Then, substituting into equation (\ref{pre_cnoidal_soln}) we obtain the periodic solution found by Korteweg and deVries:
\begin{equation}\label{cnoidal_soln_2}
u(x,t) = f(x - ct) = f_2 - (f_2 - f_3)\; \text{cn}^2\left(    \left.  \sqrt{\frac{f_1 - f_3}{2}} \left(  x - ct - x_0  \right)  \right| m \right)
\end{equation}
The period of the solution is
\[
P=     2 \int_{f_3}^{f_2}  \frac{1}{ \sqrt{2F(\hat f)}}d\hat f  = 2 K(m) {\sqrt{\frac{2}{f_1 - f_3}}},
\]
where
$$
K(m):= 2\int_{0}^{\frac{\pi}{2}}  \frac{1}{ \sqrt{1 - m \sin^2\theta}}d\theta.
$$

This is the only solution that can be computed by means of the traveling wave solution assumption. This is actually a single-phase solution of the KdV equation, with phase $\theta = x - ct - x_0$. With the work developed in \cite{Kamchatnov_Kraenkel_Umarov_1, Kamchatnov_Kraenkel, Kamchatnov_Kraenkel_Umarov_2} by Dubrovin, Novikov, Its, Faddeev. Matveev, more solutions can be developed, in particular, multi-phase solutions. We are going to follow a general scheme starting from the spectral problem (\ref{e_val_problem}) for the KdV equation and we develop multiphase solutions. 
In particular, we will recover equation (\ref{separable_eqn}), see section \ref{sect:1_soliton_soln}. 

\section{Stationary solutions of the KdV equation}\label{sect:stationary_soln_s_kdv}

Consider again the Korteweg-deVries equation (\ref{kdv}),
\begin{equation}\label{kdv_2}
q_t + 6qq_x + q_{xxx}  = 0.
\end{equation}
We know that, to the KdV equation, we have an associated linear eigenvalue problem, equations (\ref{e_val_problem}) and (\ref{L_operator}), which imply the Schr\"odinger equation
\begin{equation}\label{eval_problem_2}
y" + q y = \lambda y.
\end{equation}
Here, we consider stationary solutions of the KdV equations, $q = q(x,t_0)$ for a fixed time $t_0$, and we will recover solutions to equation (\ref{kdv_2}) from the associated Riemann surface. Next, we will consider how this Riemann surface evolves in time. We start following the ideas of Flaschka \cite{flaschka_1}.

Look for solutions of (\ref{eval_problem_2}) of the form
\begin{equation}\label{squared_ei_funct}
\phi = y^2.
\end{equation}
These functions are called {\em "squared eigenfunctions"}.  In section 7, we obtain the following differential equation for $\phi$
\begin{equation}\label{eqn_4_fi}
\phi''' + 4q\phi + 2q' \phi = 4 \lambda \phi'   .
\end{equation}


We assume that the function $\phi$ is also a polynomial function of $\lambda$,
{\em i.e.}, $\lambda$ just appears as a parameter. To set ideas, set $\phi_N = \phi_N(x; \lambda)$ a polynomial of degree $N$,  
\[
\phi_N(x; \lambda) = (4\lambda)^NF_{-1} + (4\lambda)^{N-1}F_{0}  + (4\lambda)^{N-2}F_{1}  + \dots + (4\lambda)F_{N-2}  + F_{N-1} ,
\]
thus, we prove the following lemma. 

\lmm\label{lem0}
The function $\phi_N(x; \lambda) $ solves (\ref{eqn_4_fi}) if, and only if,  $F_N$ is constant (in $x$), {\em i.e.}, 
\[
F_N = F_N(q(x), q'(x), q"(x), q''(x), \dots, q^{(2N)})  = K 
\]
is independent of $x$, 
\mml

For a detailed proof, see Corollary \ref{corollary_main}.


In the meanwhile, we will show how the result works for $N=0,1,	2$.

\section{The 0-, 1- and 2-soliton solutions for the KdV equation}\label{sect:0_1_2_soln_s_4_kdv}

\subsection{The 0--soliton solution: $N=0$}\label{sect:0_soln}
In this instance, 
\[
\phi_0(x; \lambda) = A_0(x),
\] 
is a zero-degree polynomial in $\lambda$ (here, $F_{-1} = A_0(x)$). Substitute into (\ref{eqn_4_fi}), to get
\[
 A_0''' +  4 q  A_0' + 2q'  A_0 = 4\lambda  A_0'.
\]
On the right-hand-side of this equation, we have a first degree polynomial in $\lambda$. Then:
\begin{eqnarray*}
 A_0'   & = & 0, \\
A_0''' +  4 q  A_0' + 2q'  A_0 & = &  0.
\end{eqnarray*}
Hence, $A_0(x) = A_0$ is a constant, and $2q' A_0= 0$. Assuming $A_0 \neq 0$, then we must have
$$q(x) = q_0$$
is also a constant, which is the trivial solution to the KdV equation. This is to say, this is the {\bf $0-$soliton solution}.

\subsection{The 1-soliton solution: $N=1$}\label{sect:1_soln}
Consider
\begin{equation}\label{defn_fi_1}
\phi_1(x; \lambda) = \lambda A_0 + A_1(x), 
\end{equation}
where $ A_1(x)$ has to be found. Substitute $\phi_1$ into equation (\ref{eqn_4_fi}), to get
\[
P[\phi_1] = 4\lambda \phi_1',
\]
where
\begin{equation}\label{defn_P}
P = \frac{d^3}{dx^3}  +   4 q \frac{d}{dx}   + 2q_x.
\end{equation}
We then have,
\[
A_1''' + 4qA_1' + 2q'(\lambda A_0 + A_1(x)) =  4\lambda A_1'.
\]
At $O(\lambda)$, $2A_1' = q' A_0$, setting $A_0 = 2$, we have
\begin{eqnarray}
A_1' &=& q'    \label{A_1_prime}  \\
A_1(x) &=& q(x) + q_1,  \label{A_1_equals_q} 
\end{eqnarray}
with $q_1$ is a constant of integration. Hence:
\begin{equation}\label{soln_of_fi_1}
\phi_1(x; \lambda) =   2 \lambda  + (q(x) + q_1), 
\end{equation}

At $O(\lambda^0) = O(1)$,
\begin{eqnarray*}
A_1''' + 4qA_1' + 2q'A_1 &=& 0, \\ 
q''' + 4qq' + 2q'(q + q_1) &=& 0,
\end{eqnarray*}
{\em i.e.}, 
\[
q''' + 6qq' + 2q_1q' = 0.
\]
which is the stationary KdV equation! Notice this equation can be written as:
\begin{equation}\label{1st_flux}
q'' + 3q^2 + 2q_1q = K = \text{constant},
\end{equation}
which is a constant flux for the stationary KdV equation.

\subsection{The 2-soliton solution: $N=2$}\label{sect:2_soln}

We define
\begin{equation}\label{defn_fi_2}
\phi_2(x; \lambda) = \alpha_2\lambda \phi_1+ A_2(x), 
\end{equation}
and this formula will be used as a recursion formula. Here, we require to find $A_1(x)$  and the appropriate constant   $\alpha_2$. 
To set ideas, consider $q_1 = 0$ in (\ref{soln_of_fi_1}). Use (\ref{defn_fi_1}) in (\ref{defn_fi_2}) to get:
\[
\phi_2(x; \lambda) = \alpha_2A_0\lambda^2 + \alpha_2\lambda A_1(x)+ A_2(x).
\]
Substitute this $\phi_2$ into (\ref{eqn_4_fi}), to get
\begin{eqnarray*}
P[\phi_2] &=& 4\lambda\phi_2'\\
P[\alpha_2A_0\lambda^2 + \alpha_2\lambda A_1+ A_2] &=& 4\lambda(\alpha_2A_0\lambda^2 + \alpha_2\lambda A_1+ A_2)' .
\end{eqnarray*}

Since $A_0=$constant, and by linearity of $P$, 
\[
\alpha_2\lambda^2P[A_0] +  \alpha_2\lambda P[A_1] + P[A_2]   =   4\alpha_2\lambda^2 A_1' + 4\lambda A_2'.
\]

Now, we have to determine $A_2(x)$ and $\alpha_2$. Consider terms of the order $O(\lambda^2)$ in the previous equation. We have $\alpha_2P[A_0]  =  4\alpha_2 A_1' $. Since $A_0 = 2$ is a constant, and assuming $\alpha_2 \neq 0$, we have $q_x= A_1'$, which holds identically since $ A_1(x)  = q(x) + q_1$. 

At order $O(\lambda)$, we have $\alpha_2 P[A_1]  =   4 A_2'$. Since $A_1(x)  = q(x) + q_1$, and using the definition of $P$, equation (\ref{defn_P}),
we have 
\[
\alpha_2(q''' + 4qq' + 2q'[q + q_1]) = 4A_2'.
\]
Choose $\alpha_2 = 4$. Hence $q''' + 6qq' + 2q_1q' = A_2'$ {\em i.e.}, $A_2' = (q'' + 3q^2 + 2 q_1q)'$. Hence:

\be\label{expression_4_A_2}
A_2 = q'' + 3q^2 + 2 q_1q + q_2,
\ee
where $q_2$ is another constant. This is again a constant flux for the KdV equation, as found in equation (\ref{1st_flux}).

At order $O(\lambda^0) = O(1)$, this is to say, $ P[A_2]   =  0 $, {\em i.e.}, 
\be\label{eqn_4_A_2}
A_2''' + 4qA_2' + 2q_xA_2 = 0.
\ee
Substitute (\ref{expression_4_A_2}) into  the last equation (\ref{eqn_4_A_2}), to get
\[
q^{(5)} + 10 q q''' + 20 q'q'' + 30q^2q' + 2q_1(q''' + 6qq') + 2 q_2q' = 0.
\]
The last two terms correspond to the KdV flow and the traveling wave flow. Then, setting $q_1 = q_2 = 0$, we don not loose any information.
Hence,
\be\label{5th_order_kdv}
q^{(5)} + 10 q q''' + 20 q'q'' + 30q^2q' = 0,
\ee
which is the fifth order KdV equation! Or the third order KdV in the Lax's hierarchy
\[
q^{(4)} + 5 (q')^2  +  10 qq'' + 10q^3 =  \text{constant},
\]
which is also a KdV constant flux. We then have.
\be\label{fi_2}
\phi_2(x,\lambda) = 8\lambda^2 + 4\lambda(q(x) + q_1) + (q" + 3q^2 + 2q_1q + q_2) 
\ee
If $q_1 = q_2 =0$, this function corresponds to Flaschka's solution \cite{flaschka_1}.

\section{Recovering the  solutions of the KdV}\label{sect:recovering_kdv_solns}

In this section we will study how to recover the solution to the (stationary)  KdV equation out of its corresponding Riemann Surface. We will do it for the cases $N=0,1$ and $2$. For the $N=0$ case, there is no much to do, $q(x) = q_0=$constant, and this is the trivial solution to the KdV equation.

\subsection{The 1--soliton solution: $N=1$}\label{sect:1_soliton_soln}

For the $N=1$ case, the recovering is so straightforward that the general procedure cannot be really appreciated. Setting $q_1=0$, we have,
\be\label{fi_q}
\phi_1(x;\lambda) = 2\left( \lambda + \frac{q(x)}{2} \right).
\ee
Then, the roots, in $\lambda$, for $\phi_1(x;\lambda) =0$ are $\lambda = \lambda_1(x)$ (here, we just have one), and they are functions of $x$, namely, 
\be\label{fi_lambda_1}
\phi_1(x;\lambda) = 2\left( \lambda -  \lambda_1(x)  \right).
\ee
So, if we can compute the zeros of $\phi_1(x;\lambda) =  0$, we can compute the solution to the KdV equation:,
\be\label{q_1_soliton}
q(x) = - 2 \lambda_1(x).
\ee
This may look quite trivial. But we will see how to find $\lambda_1(x)$ so that we can use equation (\ref{q_1_soliton}) to solve the KdV equation. 

Consider the equation for $\phi$, (\ref{eqn_4_fi}), and multiply it by $\phi$ itself and, after an integration by parts, we obtain
\be\label{eqn_4_fi_and_R}
\phi\phi'' - \frac{1}{2}(\phi')^2  + 2 (q - \lambda)\phi^2 =: R(\lambda),
\ee
where the right-hand-side is a polynomial in $\lambda$. Substitute (\ref{fi_q}) into (\ref{eqn_4_fi_and_R}) above, to get
\be\label{R_3}
R_3(\lambda):=-  \frac{1}{8} R(\lambda) = \lambda^3   -  \frac{1}{4}\lambda(q'' + 3q^2) +   \frac{1}{8}\left((q^3 + \frac{1}{2}(q')^2)-q(q''+3q^2) \right).
\ee
(Here, $R_3(\lambda)$ represents a polynomial of degree $3$).

Assume we substitute the initial conditions for the KdV equation, $q = q(x,0)$, into (\ref{R_3}), to have functions of $x$ only in the coefficients of the polynomial $R_3(\lambda)$. How can we recover $q(x)$ from these functions of $x$? We proceed as follows. 

From (\ref{fi_q}), it follows that $\left.  \phi\right|_{\lambda = \lambda_1}=0$. 
Hence, evaluate (\ref{eqn_4_fi_and_R}) at $\lambda = \lambda_1$ to get
\[
\left.  - \frac{1}{2}(\phi')^2\right|_{\lambda = \lambda_1(x)} = -8 R_3(\lambda_1),
\]
{\em i.e.}, 
\be\label{riemann_g_3}
\left. (\phi')^2\right|_{\lambda = \lambda_1(x)}     =   2^4 R_3(\lambda_1),
\ee
which represents a Riemann surface (an elliptic curve) of genus $g=1$ in $(\lambda, \phi')) \in \mathbb{C}^2$. On the other hand, from equation (\ref{fi_lambda_1}),
\[
\frac{d\phi}{dx}  =   - 2 \frac{d\lambda_1}{dx},
\]
so that from (\ref{riemann_g_3}) we obtain,
\[
\left.   - 2 \frac{d\lambda_1}{dx}\right|_{\lambda = \lambda_1(x)}   = 2^2 \sqrt{R_3(\lambda_1)}.
\]
{\em i.e.},
\be\label{eqn_4_fi_1}
\left.   \frac{\lambda_1'}{   \sqrt{R_3(\lambda_1)}}\right|_{\lambda = \lambda_1(x)}   =    -2,
\ee
where we have considered the positive root.  This is a separable equation for $\lambda_1(x)$. Since in this instance $R_3(\lambda)$ is a cubic equation in $\lambda$, it turns out that equation (\ref{eqn_4_fi_1}) is precisely equation (\ref{separable_eqn}). Solving for $\lambda_1(x)$, we then are able to recover $q(x)$,
\[
q(x) = -2\lambda_1(x),
\]
which is a stationary solution of the KdV equation. How to incorporate the time-dependence, will be studied later.

\subsection{The 2--soliton solution: $N=2$}\label{sect:2_soliton_soln}

In this case, we proceed as in the previous case. The idea to follow is exactly the same, although the computations are a little bit more involved. 

We have equation (\ref{fi_2}), with $q_1 = q_2 = 0$ (just to set ideas):
\be\label{fi_2_bis}
\phi_2(x,\lambda) = 8\lambda^2 + 4\lambda(q(x)) + (q" + 3q^2).
\ee
Substitute into (\ref{eqn_4_fi_and_R}) to get
\be\label{eqn_4_fi_2_and_R_5}
\phi_2\phi_2'' - \frac{1}{2}(\phi_2')^2  + 2 (q - \lambda)\phi_2^2 = -128R_5(\lambda),
\ee
where,
\[
R_5(\lambda) = \lambda^5 -\frac{1}{16}F_2\lambda^3 -\frac{1}{32}K_2\lambda-\frac{1}{128}L_2,
\]
is a fifth order polynomial in $\lambda$; and
\bas
F_2 & = & q^{(4)}+10qq''+5(q')^2+10q^3,    \\
K_2 & = &  qq^{(4)}-q'q^{(3)}+10q^2q^{(2)}+\frac{1}{2}(q^{(2)})^2+10q(q')^2+\frac{25}{2}q^4,\\
L_2 & = & q^{(4)}(q^{(2)}+3q^2)-q^{(3)}(\frac{1}{2}q^{(3)}+6qq')+q^{(2)}(8q^{(2)}q+6(q')^2+30q^3)+18q^5. 
\eas
are constants along trajectories ({\em i.e.}, solutions) of the 5$^{\text{th}}$ order KdV flow (\ref{5th_order_kdv}). Now, we can write equation (\ref{fi_2_bis}) as:
\be\label{fi_2_factored}
\phi_2(x,\lambda) = 8(\lambda  - \lambda_1(x))(\lambda  - \lambda_2(x)), 
\ee
so that 
\[
\left.  \phi_2\right|_{\lambda  = \lambda_1(x)}  = \left.  \phi_2\right|_{\lambda  = \lambda_2(x)} =0.
\]
Evaluating at $\lambda  = \lambda_j(x)$ ($j = 1,2$) equation (\ref{eqn_4_fi_2_and_R_5}), 
\be\label{eqn_4_fi_prime}
\left.  \phi_2'  \right|_{\lambda  = \lambda_j(x)}   =  16\sqrt{R_5(\lambda_j)},  \quad   j = 1,2 .
\ee
From (\ref{fi_2_factored}), it follows that
\[
\phi_2'= 8\Big( - \lambda_1'(x)(\lambda  - \lambda_2(x)) - \lambda_2'(x)(\lambda  - \lambda_1(x)) \Big), 
\]
so that
\bas
\left. \phi_2'  \right|_{\lambda  = \lambda_1(x)}    &=&    - 8 \lambda_1'(x)(\lambda_1(x)  - \lambda_2(x)) ,   \\
\left. \phi_2'  \right|_{\lambda  = \lambda_2(x)}    &=&    - 8\lambda_2'(x)(\lambda_2(x)  - \lambda_1(x)) , 
\eas
and using (\ref{eqn_4_fi_prime}), we obtain
\ba
      \frac{\lambda_1'(x)}{\sqrt{R_5(\lambda_1)}}(\lambda_1(x)  - \lambda_2(x))    &=&   -2 , \label{aux_eqn_4_lambda_1_prime}    \\
    \frac{\lambda_2'(x)}{\sqrt{R_5(\lambda_1)}}(\lambda_1(x)  - \lambda_2(x))   &=&     2.  \nonumber
\ea
Adding both equations, and dividing by $\lambda_1 - \lambda_2\neq 0$, we get 
\equ\label{eqn_diff_1}   
\frac{\lambda_1'}{\sqrt{R_5(\lambda_1)}}   +   \frac{\lambda_2'}{\sqrt{R_5(\lambda_2)}}   = 0.
\uqe
Hence,
$$
\frac{\lambda_1'}{\sqrt{R_5(\lambda_1)}}  =-  \frac{\lambda_2'}{\sqrt{R_5(\lambda_2)}} 
$$
and
\equ  \label{eqn_diff_2}
\begin{split}
\frac{\lambda_1\lambda_1'}{\sqrt{R_5(\lambda_1)}}   +   \frac{\lambda_2\lambda_2'}{\sqrt{R_5(\lambda_2)}}   &= \frac{\lambda_1\lambda_1'}{\sqrt{R_5(\lambda_1)}}   -   \frac{\lambda_2\lambda_1'}{\sqrt{R_5(\lambda_1)}} \\
&=\frac{\lambda_1'}{\sqrt{R_5(\lambda_1)}}(\lambda_1-\lambda_2)=-2
\end{split}
\uqe

Putting together equations (\ref{eqn_diff_1}-\ref{eqn_diff_2}), 
we obtain a system of differential equations for $\lambda_1$ and $\lambda_2$, namely,
\ba
\frac{\lambda_1'}{\sqrt{R_5(\lambda_1)}}   +   \frac{\lambda_2'}{\sqrt{R_5(\lambda_2)}}   = 0,  \label{1st_eqn_4_lambda_1_and_lambda_2}  \\
\frac{\lambda_1\lambda_1'}{\sqrt{R_5(\lambda_1)}}   +   \frac{\lambda_2\lambda_2'}{\sqrt{R_5(\lambda_2)}}   = -2.   \label{2nd_eqn_4_lambda_1_and_lambda_2}
\ea
If we integrate,
\ba
\int^{\lambda_1}\frac{d\lambda_1}{\sqrt{R_5(\lambda_1)}}   +  \int^{\lambda_2} \frac{d\lambda_2}{\sqrt{R_5(\lambda_2)}}   &=& C_1,  \label{1st_integral_4_lambda_1_and_lambda_2}  \\
\int^{\lambda_1}\frac{\lambda_1  d\lambda_1}{\sqrt{R_5(\lambda_1)}}   +  \int^{\lambda_2} \frac{\lambda_2   d\lambda_2}{\sqrt{R_5(\lambda_2)}}   &=& -2x + C_2,   \label{2nd_integral_4_lambda_1_and_lambda_2}
\ea
where $C_1, C_2$ are constants of integration. If we are able to integrate and invert these equations, we should be able to get
\bas
\lambda_1 &=&  \lambda_1(x),  \\
\lambda_2 &=&  \lambda_2(x).
\eas
How can we recover the solution $q(x)$ of the KdV equation out of these functions? From equation (\ref{fi_2_factored})
$$
\phi_2(x,\lambda) = 8\Big(\lambda^2  - (\lambda_1(x) + \lambda_2(x))\lambda  +  \lambda_1(x) \lambda_2(x)\Big), 
$$
and comparing to equation (\ref{fi_2_bis}), we find
\[
q(x) =  -  \frac{1}{2}(\lambda_1(x) + \lambda_2(x)),
\]
which is the solution of the KdV equation. This is the way we solve the periodic KdV equation.



\section{Generalizations in Schr\"odinger equation}\label{sec1:1}

Consider a solution $y(x,t)$ of the Schr\"odinger equation $y''+qy=\lambda y$ with
$q(x)$ and $\lambda(t)$, functions only of $x$ and $t$, respectively.

Setting  $\phi=y^2$, and taking the  derivative with respect to $x$, we obtain
$$\phi'=2yy'$$
and therefore
$$(y')^2=\frac{\left(\phi'\right)^2}{4\phi}.$$
Similarly, taking a second derivative of $\phi$, we have
\equ
\begin{split}
\phi'' &=2(y')^2+2yy''\\
&=\frac{\left(\phi'\right)^2}{2\phi}+2yy''.
\end{split} 
\uqe
Hence,
\equ\label{eq0}
yy''=\frac{1}{2}\phi''-\frac{(\phi')^2}{4\phi}
\uqe
Now, multiplying the  Schr\"odinger equation by $y$ and using (\ref{eq0}),
we get the following differential equation in terms of $\phi$:
\equ\label{eq00}
\frac{1}{2}\phi''-\frac{(\phi')^2}{4\phi}=(\lambda-q)\phi.
\uqe
A straight-forward computation shows that
\equ\label{eq000}
\begin{split}
\left(\frac{(\phi')^2}{4\phi}\right)_x &=\frac{\phi'}{\phi}\left(\frac{1}{2}\phi''-\frac{(\phi')^2}{4\phi}\right)\\
&=\frac{\phi'}{\phi}(\lambda-q)\phi=(\lambda-q)\phi' 
\end{split}
\uqe
Hence, diferentiating equation (\ref{eq00}), we obtain
\equ\label{eq0000}
\frac{1}{2}\phi'''-(\lambda-q)\phi'=\left[(\lambda-q)\phi\right]_x
\uqe
Multiplying by $2$ and re-arranging terms, we finally get the diferential equation:
\equ\label{eq1}
\phi'''+4q\phi'+2q'\phi=4\lambda \phi'.
\uqe
Define the linear operator
$$B=\frac{d^3}{dx}+4q\frac{d}{dx}+2q';$$
thus, equation (\ref{eq1}) can be re-written as 
\equ\label{eqB}
B(\phi)=4\lambda\phi'.
\uqe

We are interested in solutions to (\ref{eqB}) of polynomial form in $\lambda$:
\equ\label{equ2}
\phi_n(x,t)=A_n\lambda^n(t)+A_{n-1}\lambda^{n-1}(t)+\cdots+A_0
\uqe
where $A_n$ is constant, and $A_i(x)$ for $0\leq i <n$, functions of $x$.

Define recursively the following conservation laws $F_n(x)$:
$$F_{-1}=\frac{1}{2}$$
$$F_{n+1}=\int B(F_n) dx.$$
The first three terms are $F_{-1}=\frac{1}{2}$, $F_{0}=q$, $F_{1}=q''+3q^2$. 
Also, sometimes for convinience, we write the recursion among conservations laws in its differential
form  $F_{n+1}'=B(F_n)$.

\thm \label{thm1}
\mbox{ }

Let $\phi_n=(4\lambda)^nF_{-1}+(4\lambda)^{n-1}F_{0}+\cdots+F_{n-1}$. 
Thus $B(\phi_n)-4\lambda\phi'_n=B(F_{n-1})$.
\mht

\begin{proof}
First, observe
\equ
\begin{split}
\phi_n & :=(4\lambda)\left[(4\lambda)^{n-1}F_{-1}+(4\lambda)^{n-2}F_{0}+\cdots+F_{n-2}\right]+F_{n-1}\\
      & = 4\lambda \phi_{n-1}+F_{n-1}.
\end{split}
\uqe
The result is clear for $n=0$, because $\phi_0=F_{-1}=\frac{1}{2}$ and thus $\phi_0'=0$ and $B(\phi_0)=B(F_{-1})$.

Now, let's assume the result is valid for $n$, hence the calculation
\equ
\begin{split}
B(\phi_{n+1})-4\lambda\phi'_{n+1} &= B(4\lambda\phi_n+F_n)-4\lambda\left[4\lambda\phi_n+F_n\right]'\\
&= 4\lambda B(\phi_n)+B(F_n)-(4\lambda)^2\phi_n'-4\lambda F_n'\\
&= (4\lambda)\left[B(\phi_n)-4\lambda\phi_n'\right]+B(F_n)-4\lambda B(F_{n-1})\\
&= 4\lambda B(F_{n-1})+B(F_n)-4\lambda B(F_{n-1})=B(F_n).
\end{split}
\uqe
shows it for $n+1$.

\end{proof}

\crl \label{cor1} 
$B(\phi_n)=4\lambda\phi'_n$ if and only if 
$F_n$ is constant.
\lrc 

 \begin{proof} 
Since $F_{n}'=B(F_{n-1})$, it is an immediate consequence from Theorem \ref{thm1}.


 \end{proof}

The definition of the polynomials $\phi_n$ in Theorem \ref{thm1} may seem arbitrary, but they are enough to understand the polynomial solutions in $\lambda$ to (\ref{eq1}). The following lemma is useful to understand this claim.

\lmm\label{lem0}
Let 
$\psi_n=A_{-1}(4\lambda)^n+A_0(4\lambda)^{n-1}+\cdots+A_{n-2}(4\lambda)+A_{n-1}$
with $A_{-1}$ constant and $A_i$ functions of $x$ for $0\leq i \leq n-1$. 

If 
$B(\psi_n)=4\lambda\psi_n'$, thus $B(A_{i-1})=A_i'$ for  $0\leq i \leq n-1$ and $B(A_{n-1})=0$.
\mml

\prf
On the one hand,
using linearity of the operator $B$, we have
\equ
B(\psi_n) =B(A_{-1})(4\lambda)^n+B(A_0)(4\lambda)^{n-1}+\cdots+B(A_{n-2})(4\lambda)+B(A_{n-1}),
\uqe
on the other hand, taking derivative of $\phi_n$ and multiplying by $4\lambda$, we get
\equ
(4\lambda)\psi_n' = A_0'(4\lambda)^{n}+A_1'(4\lambda)^{n-1}+\cdots+A_{n-2}'(4\lambda)^2+A_{n-1}'(4\lambda).
\uqe
Hence, the condition $B(\psi_n)=4\lambda\psi_n'$ gives the equality of two polynomials of degree $n$ in $4\lambda$; since terms with the same power should concide, we obtain that
$B(A_{i-1})=A_i'$ and $B(A_{n-1})=0$.
\frp

\dfn \label{defpsi}
A function of the form
$$\psi_n=A_{-1}(4\lambda)^n+A_0(4\lambda)^{n-1}+\cdots+A_{n-2}(4\lambda)+A_{n-1}$$
with $A_{-1}$ constant and $A_i(x)$ functions of $x$ for $0\leq i \leq n-1$
which is solution of $B(\psi_n)=4\lambda\psi_n'$ is called an $n$-soliton.
\nfd

\thm\label{all_N_solitons}
\enu
\item Each $n$-soliton $\psi_n$ can be written as a linear combination of the basic solitons: $\phi_n$, $\phi_{n-1}$, $\ldots$, $\phi_0$.
\item The linear combination
$$\psi_n=\alpha_n\phi_n+\alpha_{n-1}\phi_{n-1}+\ldots+\alpha_0\phi_0$$ (with $\alpha_i$'s constant and $\alpha_n\neq 0$) is a $n$-soliton if and only if $$\alpha_nF_n+\alpha_{n-1}F_{n-1}+\cdots+\alpha_0F_0$$ is constant.
\une
\mht

\prf
\enu
\item By lemma (\ref{lem0}), assuming that  $\psi_n$ is written as in (\ref{def}), then $$B(A_{i-1})=A_i'$$ for  $0\leq i \leq n-1$.

Now, since $A_{-1}$ and $F_{-1}$ are constants not zero, set $K_{-1}$ such that
$$A_{-1}=K_{-1}F_{-1}.$$

Also, integrating
\equ
\begin{split}
A_0 & =\int B(A_{-1}) dx = \int B(K_{-1}F_{-1}) dx\\
&= K_{-1}\int B(F_{-1})dx = K_{-1}F_0+C_0=K_{-1}F_0+K_0F_{-1},
\end{split}
\uqe
where $C_0$ is a
constant and $K_0$ chosen such that $C_0=K_0F_{-1}$.

Proceeding by induction, we find constants $K_j$, such that
$$A_i=K_{-1}F_i+K_0F_{i-1}+\cdots+K_{i-1}F_0+K_iF_{-1}=\sum_{j=0}^{i+1}K_{j-1}F_{i-j}$$
for  $-1\leq i \leq n-1$.

Hence,

\equ
\begin{split}
\psi_n &=\sum_{i=0}^n (4\lambda)^iA_{n-1-i}\\
&= \sum_{i=0}^n (4\lambda)^i\sum_{j=0}^{n-i}K_{j-1}F_{n-1-i-j}\\
&= \sum_{j=0}^n K_{j-1}\sum_{i=0}^{n-j}(4\lambda)^iF_{n-1-i-j}\\
&=  \sum_{j=0}^n K_{j-1}\phi_{n-j}
\end{split}
\uqe
 
\item
Now, if $$\psi_n = \sum_{j=0}^n \alpha_{j}\phi_{j},$$ satisfies $B(\psi_n)=4\lambda\psi_n'$, 
using linearity
\equ
\begin{split}
B(\psi_n)-4\lambda\psi_n' 
&=\sum_{j=0}^n \alpha_{j}B(\phi_{j})-4\lambda\left[\sum_{j=0}^n \alpha_{j}\phi_{j}'\right]\\
&=\sum_{j=0}^n \alpha_{j}\left[B(\phi_{j})-4\lambda\phi_{j}'\right]\\
&=\sum_{j=0}^n \alpha_{j}B(F_{j})=0  
\end{split}
\uqe
Integrating, we prove that $\sum_{j=0}^n\alpha_{j}F_{j}$ is constant.

\une

\frp

\dfn \label{Hn_function}
Define $$\Hc(\phi):=\int\phi B(\phi)dx -2\lambda\phi^2=\phi\phi''-\frac{1}{2}(\phi')^2+2(q-\lambda)\phi^2.$$
\nfd

\pro
If $\phi\neq 0$, then $\Hc(\phi)$ is a constant if and only if $B(\phi)=4\lambda\phi'$.
\orp

\prf
If $B(\phi)=4\lambda\phi'$, then clearly $\phi B(\phi)-4\lambda\phi'\phi=0$
and integrating  $\Hc(\phi)=\int\phi B(\phi)dx -2\lambda\phi^2$ is a constant. If $\phi\neq 0$, the ``only if'' part is also true.

\frp

Set $\Hc_n:=\Hc(\phi_n)$. Observe that since $\phi_n$ is a polynomial in $\lambda$, $\Hc_n$ is also a polynomial in $\lambda$. The following result gives a explicite description of $\Hc_n(\lambda)$ in terms of the conservations laws $F_{-1}$, $F_0$, $\ldots$, $F_n$.

\thm \label{thm2}
\equ
\begin{split}
\Hc_n(\lambda) &=-\frac{(4\lambda)^{2n+1}F_{-1}^2}{2}+(4\lambda)^{n}F_{-1}F_n+(4\lambda)^{n-1}\left[F_0F_n-\int F_0'F_ndx \right]+\\
&\makebox[5mm]{ }+(4\lambda)^{n-2}\left[F_1F_n-\int F_1'F_ndx \right]+\cdots+\left[F_{n-1}F_n-\int F_{n-1}'F_n dx \right]
\end{split}
\uqe
\mht

The proof of Theorem \ref{thm2} is on Appendix \ref{apn1}.

\crl\label{corollary_main}
If $\phi_n$ is a solution of $B(\phi_n)=4\phi'_n$, then the polynomial $\Hc_n=H(\phi_n)$ has constant coefficients with respect to $x$.
\lrc

\begin{proof}
The conservation law $F_n$ is constant by Corollary \ref{cor1}.

Since $F_{-1}=\frac{1}{2}$, the coeficients of $\lambda^{2n+1}$ and
$\lambda^n$ are clearly constant.

For the other coeficients with $0\leq j\leq n-1$, since $F_n$ is constant, we have that
$$F_{j}F_n-\int F_{j}'F_n dx = F_n\left[F_{j}-\int F_{j}'dx\right]$$
is also constant, by the Fundamental Theorem of Calculus.

\end{proof}

\section{Geometry of Solutions}\label{sec2}

In this section, the explicit expresions of $\phi_n$ and $\Hc_n$ as polynomials of $\lambda$ obtained in section \ref{sec1:1} will be used to analize the geometry in the solutions of equation (\ref{eqB}). 

The solutions of equation (\ref{eqB}) by Theorem \ref{thm1} have the explicit form

$$\phi_n=\frac{4^n}{2}\lambda^n+4^{n-1}F_{0}\lambda^{n-1}+\cdots+4F_{n-1}\lambda+F_{n-1}$$
provided that $F_n$ is constant.

But, the conservation laws $F_0$,$F_1$, $\ldots$, $F_{n-1}$ are not constant, they depend on
the position $x$ and the time $t$. Therefore, $\phi_n(x,t)$ is a polynomial which their coeficients
vary depending on $(x,t)$. Nevertheless, we can factorize $\phi_n$ over $\C$ for fixed
values of $x$ and $t$.
\equ\label{facto}
\phi_n=\frac{4^n}{2}\prod_{i=1}^{n}\left[\lambda-\lambda_i(x,t)\right].
\uqe
Taking the derivative of the above expression with respect to $x$, we obtain
\equ
\phi'_n=-\frac{4^n}{2}\sum_{j=1}^n\lambda'_j\prod_{i\neq j}\left[\lambda-\lambda_i(x,t)\right],
\uqe
evaluating for $\lambda=\lambda_k(x,t)$, we finally get
\equ\label{eqg1}
\phi'_n\mid_{\lambda=\lambda_k}=-\frac{4^n}{2}\lambda'_k\prod_{i\neq k}\left[\lambda_k(x,t)-\lambda_i(x,t)\right].
\uqe

Now, using that
 $$\Hc_n(\lambda)=\Hc(\phi_n)=\phi_n\phi_n''-\frac{1}{2}(\phi'_n)^2+2(q-\lambda)\phi_n^2$$
and that $\phi_n\mid_{\lambda=\lambda_k(x,t)}=0$, we obtain
\equ\label{eqg2}
\Hc_n(\lambda_k)=-\frac{1}{2}(\phi'_n\mid_{\lambda=\lambda_k})^2
\uqe

Now, since $\Hc_n$ is a polynomial of degree $2n+1$ with constant coefficients with respect to the variable $x$, the equation
\equ
\Hc_n(X)=-\frac{1}{2}Y^2
\uqe
is of an hyperelliptic curve $\mathcal{H}_n$ of genus $n$ and each $P_k=(\lambda_k,\phi_n'\mid_{\lambda=\lambda_k})$ with $1\leq k\leq n$ represents a point on it. Hence, if we vary $x$, but fixed $t$, a soliton solution defines $n$ real curves on $\mathcal{H}_n$.
These curves on the tangent space define a system of linear differentials. 

In fact, combining equations (\ref{eqg1}) and (\ref{eqg2}), we have that

\equ\label{eqg3}
\begin{split}
\sqrt{-2\Hc_n(\lambda_k)} &=\phi'_n\mid_{\lambda=\lambda_k}\\
&= -\frac{4^n}{2}\lambda'_k\prod_{i\neq k}\left[\lambda_k(x,t)-\lambda_i(x,t)\right].
\end{split}
\uqe

But, since $F_{-1}^2=\frac{1}{4}$, the leading coefficient of $\Hc_n$ is $$-\frac{4^{2n}}{2}.$$

Hence, setting $$\mathcal{P}_n=\frac{-2\Hc_n}{4^{2n}},$$ we obtain
\equ\label{eqg4}
\frac{\lambda_k'}{\sqrt{\mathcal{P}_n(\lambda_k)}}=\frac{-2}{\prod_{i\neq k}\left[\lambda_k(x,t)-\lambda_i(x,t)\right]}
\uqe

Now, the differentials 
\equ\label{diff}
\omega_\mu=\frac{X^{\mu-1}dX}{\sqrt{\mathcal{P}_n(X)}}\mbox{ , $1\leq\mu\leq n$}
\uqe
form a basis of the differential space $\Omega^1(\mathcal{H}_n)$.

Evaluating those differentials in the points $P_k$, we obtain
\equ
\begin{split}
\omega_\mu(P_k)=\omega_\mu(\lambda_k) &= \frac{\lambda_k^{\mu-1}\lambda_k'dx}{\sqrt{\mathcal{P}_n(X)}}\\
&=\frac{-2\lambda_k^{\mu-1}}{\prod_{i\neq k}\left[\lambda_k(x,t)-\lambda_i(x,t)\right]} 
\end{split}
\uqe

Adding up over all points $P_k$ and using the main proposition \ref{mainapn2} in Appendix \ref{apn2}, we get
\equ\label{maindiff}
\sum_{k=1}^{n}\omega_\mu(P_k)=
\begin{cases}
0 & \mbox{ if }1\leq \mu< n \cr
-2 & \mbox{ if }\mu=n
\end{cases}
\uqe

Now, the Abel map from a Riemann surface $\mathcal{H}_n$ to its Jacobian $J(\mathcal{H}_n)$ (chosing a fixed point $P_0$) is defined by
\equ
A_{P_0}(Q):=(A_{P_0}^{(1)}(Q),A_{P_0}^{(2)}(Q),\ldots,A_{P_0}^{(n)}(Q))
\uqe
where
$$A_{P_0}^{(k)}(Q):=\int_{P_0}^{Q}\omega_k.$$

If $S_n=\bigwedge^n\mathcal{H}_n$ is the $n$-symetric power of $\mathcal{H}_n$, the Abel is the extended map
$$\tilde{A}_{P_0}:S_n\rightarrow J(\mathcal{H}_n)$$
given by the addition formula
\equ
\tilde{A}_{P_0}(Q_1,Q_2,\ldots,Q_n)=\sum_{k=1}^nA_{P_0}(Q_k)
\uqe

Hence,  for a soliton $\phi_n(x,t)$,  using equation (\ref{maindiff}), taking $$P_k(x,t)=\left(\lambda_k(x,t),\phi_n'\mid_{\lambda=\lambda_k}(x,t)\right),$$ and considering $P_0=(\lambda_0,\mu_0)$ any other point in $\mathcal{H}_n$, we obtain that
\equ
\tilde{A}_{P_0}\left(P_1(x,t),P_2(x,t),\ldots,P_n(x,t)\right)=\left(C_1(t),C_2(t),\ldots,C_{n-1}(t),C_n(t)-2x\right)
\uqe
where the functions $C_i(t)$ are constant with respect to $x$.

\section{Time Evolution}\label{sec3}

The Schr\"odinger differential operator
$$\mathcal{L}=\frac{d^2}{dx}+q$$
has an associatted operator
$$
P=-4\frac{d^3}{dx}-6q\frac{d}{dx}-3q',
$$
so that they form the  Lax pair of the Korteweg-deVries equation. 

Define the commutator operator $[\mathcal{L},P]:=\mathcal{L}P-P\mathcal{L}$. The following result describes 
the relationship among $\mathcal{L}$ and $P$.

\lmm
For $\mathcal{L}$ and $P$ as above, we have:
\equ
[\mathcal{L},P]=q_{xxx}+6qq_{x}
\uqe
\mml
\prf
To simplify our notation, set $$D:=\frac{d}{dx}.$$ If $f$ is a differentaible function over some domain $\Omega$, $D$ acts according
to the formalism:
$$Df=f_x+fD.$$
In particular, it can be proven by induction that:
$$D^nf=\sum_{i=0}^n \binom{n}{i} f^{(n-i)}D^i.$$
Thus, using this notation, we have that
\equ
\begin{split}
[\mathcal{L},P] &=[D^2+q,-4D^3-6qD-3q_x]\\
&=-\{6[D^2,qD]+3[D^2,q_x]+4[q,D^3]+6[q,qD]\},
\end{split}
\uqe
since $[D^2,D^3]=0$ and $[q,q_x]=0$.
Now,
\equ\label{conmutator1}
\begin{split}
6[D^2,qD] &=6\{D^2(qD)-(qD)D^2\}=6\{(D^2q)D-qD^3\}\\
&= 6\{(q_{xx}+2q_xD+qD^2)D-qD^3\}=6q_{xx}D+12q_xD^2
\end{split}
\uqe
\equ\label{conmutator2}
\begin{split}
3[D^2,q_x] &=3\{D^2(q_x)-(q_x)D^2\}\\
&= 3\{(q_{xxx}+2q_{xx}D+q_xD^2)-q_xD^2\}=3q_{xxx}+6q_{xx}D
\end{split}
\uqe
\equ\label{conmutator3}
\begin{split}
4[q,D^3] &=4\{qD^3-D^3q\}= 4\{qD^3-(q_{xxx}+3q_{xx}D+3q_xD^2+qD^3)\}\\
&=-4q_{xxx}-12q_{xx}D-12q_xD^2
\end{split}
\uqe
\equ\label{conmutator4}
\begin{split}
6[q,qD] &=6\{q^2D-qDq\}= 6\{q^2D-q(q_{x}+qD)\}\\
&=-6qq_{x}
\end{split}
\uqe
Adding up the equations (\ref{conmutator1}-\ref{conmutator4}), we obtain the result:
$$[\mathcal{L},P]=q_{xxx}+6qq_{x}$$

\frp

The Schr\"odinger equation, at the beggining of section \ref{sec1:1}, is the spectral or eigenfunction problem
\equ\label{spec}
\mathcal{L}y=\lambda y,
\uqe
while the $P$ operators gives the time evolution condition
\equ\label{timev}
P y =y_t=\frac{ dy }{dt}.
\uqe

As in section \ref{sec1:1}, if we consider $\phi=y^2$, we can obtain a differential equation
for $\phi_t$. In fact,
\equ\label{timev1}
\begin{split}
\phi_t &= 2yy_t=2yPy\\
&= -8yy'''-12qyy'-6q'y^2\\
&= -8yy'''-6q\phi'-6q'\phi\\
&= -8yy'''-6(q\phi)'
\end{split}
\uqe
But, using that $yy'''=(yy'')'-y'y''$ and that $y''=(\lambda-q)y$ (Schr\"odinger equation),
we obtain that
\equ\label{timev2}
\begin{split}
yy''' &=\left((\lambda-q)y^2\right)'-(\lambda-q)yy'\\
&==\left((\lambda-q)\phi\right)'-\frac{(\lambda-q)\phi'}{2}\\
&=-q'\phi+\frac{(\lambda-q)\phi'}{2}
\end{split}
\uqe
Hence,
\equ\label{timev3}
\begin{split}
\phi_t &= -8\left[-q'\phi+\frac{(\lambda-q)\phi'}{2}\right]-6(q\phi)'\\
&= 2q'\phi-2q\phi'-4\lambda\phi'\\
\end{split}
\uqe 
But, using that
$4\lambda\phi'=B(\phi)=\phi'''+4q\phi'+2q'\phi$, we obtain the linearized Kdv equation
\equ\label{timev4}
\phi_t=-\phi'''-6q\phi'
\uqe

Now, let's assume that $\phi$ is the $n$-soliton $\phi_n$.
Hence, evaluating $\phi_t$ in $\lambda=\lambda_k$, similarly as in section \ref{sec2} with $\phi'$
and using (\ref{timev3}), we obtain that
\equ\label{timev5}
\phi_t\mid_{\lambda=\lambda_k}=-2(q+2\lambda)\phi'\mid_{\lambda=\lambda_k}
\uqe

A similar computation as the one in equation (\ref{eqg1}) gives
\equ\label{eqgt}
\phi_t\mid_{\lambda=\lambda_k}=-\frac{4^n}{2}\frac{\partial \lambda_k}{\partial t}\prod_{i\neq k}\left[\lambda_k(x,t)-\lambda_i(x,t)\right].
\uqe
Hence, equation (\ref{timev5}) simplifies to the following relation of the partial derivatives of $\lambda_k$
\equ\label{timev6}
\frac{\partial \lambda_k}{\partial t}=-2(q+2\lambda)\frac{\partial \lambda_k}{\partial x}
\uqe

\part{The NLS equation}\label{part2}

\section{The Lax pair for the NLS equation. Analogies with the KdV equation.   }\label{sect:analogies_and_Lax_pair}

We know that the Nonlinear Schr\"odinger (NLS) equation
\be\label{NLS}
i\frac{\partial q}{\partial t} +  \frac{1}{2}\frac{\partial^2 q}{\partial x^2}  + \sigma |q|^2 q=   0,
\ee
(where $\sigma = \pm 1$ is the focusing/defocusing parameter) is a completely integrable system by means of the inverse scattering transform \cite{zajarov, newell, drazin}.
In this context, the NLS admits a Lax's pair in matrix form
\ba
\mathcal{L}    & = &  i \lambda D + N,  \label{L_lax_pair} \\
\mathcal{P}     & = &   i \lambda \mathcal{L} + \frac{1}{2} D(N_x - N^2) \label{P_lax_pair}
\ea
where $D$ and $N$ are two $2\times2$ matrices over $\C$, defined by
$$  D =   \left(\begin{matrix} -1 & 0 \\ 0 & 1  \end{matrix}\right)  \qquad \text{and} \qquad   N  =   \left(\begin{matrix} 0 & q \\ -\sigma \bar{q} & 1  \end{matrix}\right),$$
with $q = q(x,t)$ is solution to the NLS equation  and $\bar{q}$ denotes the complex conjugate of $q$. 

The Lax's pair, $\mathcal{L}  , \mathcal{P}  $, defines an overdetermined system of ordinary differential equations,
\ba
 \frac{\partial \Psi}{\partial x}   & = &   \mathcal{L}  \Psi,  \quad  (\text{spectral problem})  \label{ZS_system} \\
i  \frac{\partial \Psi}{\partial t}  & = &  \mathcal{P}    \Psi,  \quad  (\text{t-flow})     \label{ZS_evolution_time}
\ea
whose compatibility condition, $(\Psi_x)_t = (\Psi_t)_x$ (cross-differentiation holds) becomes the NLS equation, equation (\ref{NLS}) above.


\section{The  scalar spectral problem for the squared eigenfunctions}\label{sect:scalar_spectral_squared_eigenfunctions}

Computations made from Sections \ref{sect:stationary_soln_s_kdv} to  \ref{sec1:1} can be generalized under the hypothesis of having a linear multiplicative operator
$\hat{\mathcal{L}}$ satisfying the differential equation:

\equ
\frac{\partial^2y}{\partial x^2}=\hat{\mathcal{L}}y
\uqe

Multiplying by $y$ and using that $\hat{\mathcal{L}}$ is a multiplicative operator, we get
\equ
y\frac{\partial^2y}{\partial x^2}=\hat{\mathcal{L}}y^2
\uqe

Hence, setting $\phi=y^2$ and using equation  (\ref{eq0}), we obtain
\equ\label{rec:1}
\frac{1}{2}\phi''-\frac{(\phi')^2}{4\phi}=\hat{\mathcal{L}}\phi
\uqe

Taking derivatives in both sides, we obtain:
\equ\label{rec:1}
\frac{1}{2}\phi'''-\left(\frac{(\phi')^2}{4\phi}\right)_x=\left[\hat{\mathcal{L}}\phi\right]_x
\uqe
 
But, now repeating computations made in (\ref{eq000}) and (\ref{eq0000}), we obtain
\equ\label{rec:2}
\left(\frac{(\phi')^2}{4\phi}\right)_x =\frac{\phi'}{\phi}\left(\frac{1}{2}\phi''-\frac{(\phi')^2}{4\phi}\right)=\frac{\phi'}{\phi}\hat{\mathcal{L}}\phi=\hat{\mathcal{L}}\phi' 
\uqe
and 
\equ\label{rec:3}
\frac{1}{2}\phi'''-\hat{\mathcal{L}}\phi'=\left[\hat{\mathcal{L}}\phi\right]_x
\uqe

After multiplication by $2$, we obtain the third order differential equation
\equ\label{rec:4}
\phi'''-2\hat{\mathcal{L}}_x\phi-4\hat{\mathcal{L}}\phi'=0.
\uqe
In the case of KdV, setting $\hat{\mathcal{L}}=\lambda-q$, the equation (\ref{eq1}) is obtained, which is the basis for the KdV's recursion formul{\ae}.

In what follows, we will use a similar approach as the one used in Sections \ref{sect:stationary_soln_s_kdv}-\ref{sec1:1},   to obtain a recursion formul{\ae} for the NLS equation: 
\equ\label{nls}
iq_t+\frac{1}{2}q_{xx}+\sigma\|q\|^2q_x=0.
\uqe

The basic assumption is, 
as already mentioned, the fact that we have a scalar multiplicative operator $\hat{\mathcal{L}}$. We will make use of the work of Kamchatnov, Kraenkel and Umarov \cite{Kamchatnov_Kraenkel_Umarov_1, Kamchatnov_Kraenkel, Kamchatnov_Kraenkel_Umarov_2}, where the authors find the scalar multiplicative operator $\hat{\mathcal{L}}$ associated to the NLS equation. 
 
The scalar multiplicative operator is:
\equ
\begin{split}
\hat{\mathcal{L}} &=-\left(\lambda-\frac{iq_x}{2q}\right)^2-\sigma\|q\|^2-\left(\frac{q_x}{2q}\right)_x\\
&= -\lambda^2+E\lambda+F
\end{split},
\uqe
with 
\equ\label{Eexpresion}
E=\frac{iq_x}{q} 
\uqe
and 
\equ\label{Fexpresion}
F=-\frac{1}{4}E^2-\sigma\|q\|^2+\frac{i}{2}E', 
\uqe
which is a polynomial expresion of degree two in $\lambda$.

To simplify our computations,  define the following bilinear operator:
\equ\label{bilinear:operator}
\begin{split}
\langle \psi, \phi \rangle &= (\psi\phi)_x+\psi\phi_x\\
&= \psi_x\phi+2\psi\phi_x
\end{split}.
\uqe
Hence, equation (\ref{rec:4}) can be written as:
\equ\label{recursion:1}
\phi'''-2\langle\hat{\mathcal{L}} , \phi \rangle=0.
\uqe

Using linearity, and that $\lambda$ is constant with respect to $x$, (\ref{rec:4}) becomes:
\equ\label{recursion:2}
\phi'''-2\left[-\lambda^2\langle 1, \phi \rangle+\lambda\langle E, \phi \rangle+\langle F, \phi \rangle\right]=0.
\uqe
But, using that $\langle 1, \phi \rangle=2\phi'$, we obtain:
\equ\label{recursion:3}
\phi'''+4\lambda^2\phi'-2\left[\lambda\langle E, \phi \rangle+\langle F, \phi \rangle\right]=0.
\uqe

This equation is a scalar spectral problem for the squared-eigenfunction $\phi$.

\section{Recursion Formul\ae}\label{section:recursion}

Now, as we did in the case of the KdV equation, we will assume that $\phi$ is polynomial in $\lambda$, and that $\lambda$ is constant with respecto to $x$.
Set
$$\phi=A_0\lambda^n+A_1\lambda^{n-1}+\cdots+A_{n-1}\lambda+A_n,$$
without lost of generality, we assume $A_0$ is constant. (In fact, it is an easy exercise, to verify that if $\phi$ satisfies the equation \ref{recursion:3}, then $A_0$ is constant.)

Hence, taking derivatives, we have
$$\lambda^2\phi'=A_1'\lambda^{n+1}+A_2'\lambda^{n}+\cdots+A_{n-1}'\lambda^3+A_n'\lambda^2$$
and
$$\phi'''=A_1'''\lambda^{n-1}+A_2'''\lambda^{n-2}+\cdots+A_{n-1}'''\lambda+A_n'''.$$

Now, using that the operator $\langle \_ , \_ \rangle$ is bilinear, we obtain
$$\lambda\langle E, \phi \rangle=\lambda^{n+1}\langle E, A_0\rangle+\lambda^{n}\langle E, A_1\rangle+\cdots+\lambda^2\langle E, A_{n-1}\rangle+\lambda\langle E,A_n \rangle$$
$$\langle F, \phi \rangle=\lambda^{n}\langle F, A_0\rangle+\lambda^{n-1}\langle F, A_1\rangle+\cdots+\lambda\langle F, A_{n-1}\rangle+\langle F,A_n \rangle$$

Next, we will compare orders in equation (\ref{recursion:3}).

Comparing in order $\lambda^{n+1}$:
\equ\label{recursionls:0} 
2A_1'=\langle E, A_0\rangle
\uqe
and thus,
\equ\label{recursionls:1} 
A_1=\frac{\int \langle E, A_0\rangle dx}{2}
\uqe

Comparing in order $\lambda^{n}$:
\equ\label{recursionls:2}  
2A_2'=\langle E, A_1\rangle+\langle F, A_0\rangle
\uqe
and hence,
\equ \label{recursionls:3} 
A_2=\frac{1}{2}\int\left[\langle E, A_1\rangle+\langle F, A_0\rangle\right] dx
\uqe

Comparing in order $\lambda^{n-1}$, we obtain:
\equ \label{recursionls:4} 
A_3=\frac{1}{2}\int\left[\langle E, A_2\rangle+\langle F, A_1\rangle\right] dx-\frac{1}{4}A_1''
\uqe

In general, comparing in order $\lambda^{n-k}$ for $1\leq k \leq n-2$, we obtain
\equ 
A_{k+2}=\frac{1}{2}\int\left[\langle E, A_{k+1}\rangle+\langle F, A_k\rangle\right] dx-\frac{1}{4}A_k''
\uqe

Notice that the last comparison for $k=n-2$ occurs in order $\lambda^2$.

If we shift the indexes by $-2$, we obtain the recursion formula:
\equ \label{recursionls:5} 
A_{j}=\frac{1}{2}\int\left[\langle E, A_{j-1}\rangle+\langle F, A_{j-2}\rangle\right] dx-\frac{1}{4}A_{j-2}''
\uqe
If we define $A_{-1}=0$, this recursion formula is valid for $j=1,\ldots,n$.

Finally, comparing terms in $\lambda$ and $1$, we obtain two conditions for $\phi$ to be a solution of the NLS equation:

\enu
\item[]
\equ\label{condition:A}
A_{n-1}''=2\int\left[\langle E, A_n\rangle+\langle F, A_{n-1}\rangle\right] dx\makebox[1.5in]{ (Condition A)}
\uqe

\item[]
\equ\label{condition:B}
A_{n}''=2\int\langle F, A_n\rangle dx\makebox[1.5in]{ (Condition B)}
\uqe 
\une

\section{NLS N-solitons.}

As in the KDV case, we define some basic polynomials using the recursion formula (\ref{recursionls:5})  to characterize
all polynomial solutions to equation (\ref{recursion:3}).

Set $$\phi_0:=A_0=2.$$
Thus, 
\equ
\begin{split} 
A_1 &=\frac{\int \langle E, 2\rangle dx}{2} \\
 &= \frac{1}{2}\int 2 E_x dx=E+C
\end{split}
\uqe
with $C$ a constant.
If we set $C=0$, we obtain
$$\phi_1:=2\lambda+E.$$

Continuing with the recursion, we obtain
\equ 
\begin{split}
A_2 &=\frac{1}{2}\int\left[\langle E, E\rangle+\langle F, 2\rangle\right]\\
&=\frac{1}{2}\int\left[3EE'+2F'\right]\\
&=\frac{3}{4}E^2+F+C
\end{split}
\uqe

Taking again $C=0$, we get
$$\phi_2:=2\lambda^2+E\lambda+\frac{3}{4}E^2+F$$

In general, we define the basic NLS N-soliton $\phi_n$ by using the recursion 
formula (\ref{recursionls:5}), setting $A_0=2$ and taking all constants of integration equal to zero.

Next terms are

$$\phi_3:=2\lambda^3+E\lambda^2+(\frac{3}{4}E^2+F)\lambda+\frac{5}{8}E^3+\frac{3}{2}FE-\frac{1}{4}E''$$

\begin{align}
\phi_4:= & 2\lambda^4+E\lambda^3+\left(\frac{3}{4}E^2+F\right)\lambda^2+\left(\frac{5}{8}E^3+\frac{3}{2}FE-\frac{1}{4}E''\right)\lambda \cr
& +\left(\frac{35}{64}E^4+\frac{15}{8}E^2F+\frac{3}{4}F^2-\frac{5}{16}\left(E'\right)^2-\frac{5}{8}EE''-\frac{1}{4}F''\right) 
\end{align}

Now, suppose that

$$\psi_n(\lambda)=B_{0}\lambda^n+B_1\lambda^{n-1}+\cdots+B_{n-1}\lambda+B_{n}$$
(with $B_0$ a constant and $B_i$ function of $x$ and $t$) is a solution
of equation \ref{recursion:3}. 
We will call such a function a NLS N-soliton. The following theorem characterize such functions. It is an analogue of theorem \ref{all_N_solitons} for the NLS-equation.

\thm
\enu
\item Any $\psi_N$ can be written as a linear combination of the basic NLS-N solitons: $\phi_N$, $\phi_{N-1}$,$\ldots$,$\phi_0$.
\item Moreover, if 
$$\psi_n=\alpha_0\phi_n+\alpha_{1}\phi_{n-1}+\ldots+\alpha_n\phi_0$$
then 
$$B_i=\sum_{j=0}^i \alpha_j A_{i-j}$$
where the $B_i$'s are the coefficients of $\psi_N$ and the $A_i$'s are the coefficients of the $\phi_k$'s. Hence, the conditions (A) and (B) can be written as
$$\sum_{k=0}^n \alpha_k \mathcal{A}_{n-k}=0 \mbox{ (Condition A)}$$
and
$$\sum_{k=0}^n \alpha_k \mathcal{B}_{n-k}=0 \mbox{ (Condition B)}$$
where
$$\mathcal{A}_i=A_{i-1}''-2\int\left[\langle E, A_i\rangle+\langle F, A_{i-1}\rangle\right] dx$$
and
$$\mathcal{B}_i=A_{n}''-2\int\langle F, A_n\rangle dx$$
for $i=0,1,\ldots,n$. Assuming $A_{-1}=0$.
\une
\mht

\section{The 0-, 1- and 2- soliton solutions for the NLS equation}\label{sect:0_1_2_soliton_solns}

In this section, we explain how the recursion formulas obtained in the previous section are used. We will develope the cases when
$n=0$, $1$, $2$.

\subsection{0-soliton for NLS}

We will assume $\phi=A$ is constant and $\hat{\mathcal{L}} =-\lambda^2+E\lambda+F$. 
Thus, equation \ref{recursionls:0} gives
$$ 
\langle E, A_0\rangle=A_0E'=0
$$
Hence, $$E=\frac{iq_x}{q}=i\left(\ln{q}\right)'=k$$ is constant.
So, $q$ satisfies the linear equation
$$q_x=-ikq.$$

Thus, it follows that $$q=Ce^{-ikx},$$ with $C$ and $k$ constants.

Now, equation (\ref{recursionls:2}) gives
$$ 
\langle F, A_0\rangle=A_0F'=0
$$
Hence, $F$ is constant. But,
\equ
\begin{split}
F &=-\frac{1}{4}E^2-\sigma\|q\|^2+\frac{i}{2}E'\\
&=-\frac{1}{4}k^2-\sigma\|q\|^2.
\end{split}
\uqe
Hence, $\|q\|^2$ is constant, since $F$, $E=k$ and $\sigma$ are constants.

Thus, we can conclude that $q=Ce^{-ikx}$ with $k\in\R$.

\subsection{1-soliton for NLS}

For $n=1$, equation (\ref{recursionls:2}) gives Condition (A):
\equ\label{condition:A1}  
\langle E, A_1\rangle+\langle F, A_0\rangle =E_xA_1+2EA_1'+F_xA_0=0
\uqe

And, Condition (B) in its differential form is
\equ\label{condition:B1}
A_{1}'''-2\langle F, A_1\rangle=A_1'''-2F_xA_1-4FA'_1=0
\uqe 

But, from the recursion formula  (\ref{recursionls:1})
\equ\label{recursion:R1}
\begin{split} 
A_1 &=\frac{\int \langle E, A_0\rangle dx}{2} \\
 &= \frac{1}{2}\int A_0E_x dx
\end{split}
\uqe

Hence, taking $A_0=2$ and integrating, we obtain:
$$A_1=E$$

Hence, Condition (A) in (\ref{condition:A1}) becomes
\equ\label{condition:A1:1}
3EE_x+2F_x=0.
\uqe
Integrating with respect to $x$, we obtain
\equ\label{condition:A1:1}
\frac{3}{2}E^2+2F=C,
\uqe
where $C$ is a constant. 
Now, using (\ref{Fexpresion}, \ref{E2:definition}, \ref{E:first:derivative}), we can compute
\equ\label{2Fexpresion}
\begin{split}
2F &= -\frac{1}{2}E^2-2\sigma\|q\|^2+iE'\\
&=-\frac{1}{2}E^2-2\sigma\|q\|^2+i\left[E_{(2)}+iE^2\right]\\
&=-\frac{3}{2}E^2-2\sigma\|q\|^2+iE_{(2)}
\end{split}
\uqe

Hence, setting $-2\omega=C$, we get
\equ\label{condition:A1:2}
\frac{3}{2}E^2+2F=-2\sigma\|q\|^2-\frac{q_{xx}}{q}=-2\omega.
\uqe

Multiplying by $-\frac{q}{2}$, we obtain the Stacionary Non-linear Shrodinger equation:
$$\sigma\|q\|^2q+\frac{1}{2}q_{xx}=\omega q .$$ 

Similarly, Condition B in  \ref{condition:B1} becomes 
\equ\label{condition:B1:1}
E'''-2F_xE-4FE'=0
\uqe
using that $F'=-\frac{3}{2}EE'$ and $F=-\omega-\frac{3}{4}E^2$, we obtain
\equ\label{condition:B1:2}
E'''+6E^2E'+4\omega E'=0
\uqe
Integrating:
\equ\label{condition:B1:3}
E''+2E^3+4\omega E=C_1
\uqe
Multiplying by $E'$ and integrating again:
\equ\label{condition:B1:4}
\frac{1}{2}(E')^2+\frac{1}{2}E^4+2\omega E^2=C_1E+C_2
\uqe
Hence, we can parametrize the solutions as pair of points $(E,E')$
in the curve of genus $1$
\equ\label{nlscurve1}
y^2=-x^4+4\omega x^2+C_1x+C_2
\uqe
In fact, the above equation defines a family of elliptic curves, since $C_1$
and $C_2$ are constants. 

\subsection{2-soliton for NLS}

Using the recursion formulas, we can take $A_0=2$ and $A_1=E$, thus, by (\ref{recursionls:3})  
\equ \label{recursion:R2} 
\begin{split}
A_2 &=\frac{1}{2}\int\left[\langle E, E\rangle+\langle F, 2\rangle\right]\\
&=\frac{1}{2}\int\left[3EE'+2F'\right]\\
&=\frac{3}{4}E^2+F=-\omega
\end{split}
\uqe
But, this $\omega$ is not constant, it is variable with respect to $x$.

Now, Condition (A) is
\equ\label{condition:A2}
\begin{split}
E'' &=2\int\left[\langle E, -\omega\rangle+\langle F, E\rangle\right]\\
&=2\int\left[\langle E, \frac{3}{4}E^2+F\rangle+\langle F, E\rangle\right]\\
&=2\int\left[\frac{3}{4}\langle E,E^2\rangle+\langle E, F\rangle+\langle F, E\rangle\right]
\end{split}
\uqe
But,
\equ
\int\left[ \langle E, F\rangle+\langle F, E\rangle\right]=\int 3(FE)_x=3FE
\uqe
and
\equ
\begin{split}
\int \langle E, E^2\rangle &=\int\left[E_xE^2+2E(2E\cdot E_x)\right]\\
&=5\int \left(E^2E_x\right)=\frac{5}{3}E^3
\end{split}
\uqe

Hence, Condition (A) simplifies to
\equ\label{condition:A2:1}
E''=\frac{5}{2}E^3+6FE
\uqe

Now, using the explicite expression (\ref{Fexpresion}) for $F$ in terms of $E$, we obtain:
\equ\label{condition:A2:2}
E''=E^3-6\sigma\|q\|^2E+3iE'E
\uqe

But, using that $E'=E_{(2)}+iE^2$, equation (\ref{E:first:derivative}),
we obtain
\equ\label{condition:A2:3}
E''=-6\sigma\|q\|^2E+3iE_{(2)}E-2iE^3
\uqe
Comparing with expresion for $E''$ in (\ref{E:second:derivative}),
we finally get
\equ\label{condition:A2:4}
E_{(3)}=-6\sigma\|q\|^2E,
\uqe
which can be express in terms of $q$ and its derivatives, after multiplication by $q$ and
considering a constant of integration $C$ by the equation:
\equ\label{condition:A2:5}
iq_{xxx}+6\sigma\|q\|^2q_x=Cq.
\uqe

\section{Summary: NLS}

We follow a similar aproach as the one used in the KDV equation.

For  the NLS equation $iq_t+\frac{1}{2}q_{xx}+\sigma\|q\|^2q_x=0$, the operator

$$\hat{\mathcal{L}}=-\left(\lambda-\frac{iq_x}{2q}\right)^2-\sigma\|q\|^2-\left(\frac{q_x}{2q}\right)_x$$

is the scalar multiplicative operator for the equation:
$$\mbox{ }\mbox{ }y''=\hat{\mathcal{L}}y.$$



If we write $\hat{\mathcal{L}}$ as a polynomial in $\lambda$, we obtain
$$\hat{\mathcal{L}}=-\lambda^2+E\lambda+F,$$
where $E=\frac{iq_x}{q}$ and $F=-\frac{1}{4}E^2-\sigma\|q\|^2+\frac{1}{2}E_x$.

Using again the geometric condition $\phi=y^2$, assuming $\phi$ is of polynomial in $\lambda$, we solve the equation: 
$$\mbox{ (G) }\mbox{ }\mbox{ }\phi'''-4 \hat{\mathcal{L}}\phi'-2\hat{\mathcal{L}}'\phi=0$$
depending on the parameters $E$ and $F$.

(Condition (G) follows from equation (A) and the geometric condition)

Solving as in the case of KDV for a polynomial function in $\lambda$ of degree $n$:
$$\phi_n=A_0\lambda^n+A_1\lambda^{n-1}+\cdots+A_{n-1}\lambda+A_n,$$
 we obtain for the solutions of (G):

\enu
\item A recursion formulae:
$$A_j=\frac{1}{2}\left[EA_{j-1}+FA_{j-2}\right]+\frac{1}{2}\int\left[EA'_{j-1}+FA'_{j-2}\right]dx$$
with the initial conditions $A_{-1}=0$ and $A_0=2$.
The choice of the constant $A_0=2$ makes computations easier, but one can consider any other constant $A_0$.
Besides this recusrsio

\item Two Extra Conditions in the coefficients:

Defining the bilinear operator 
$$\langle \psi, \phi \rangle = (\psi\phi)_x+\psi\phi_x,$$
the conditions are:

\emph{Condition A}
$$A'''_{n-1}-2\langle E, A_n\rangle -2\langle F, A_{n-1}\rangle = 0$$

\emph{Condition B}

$$ A'''_n-2\langle F, A_n \rangle=0$$

Conditions (A) and (B) can also be expressed in an integral form, useful in some computations. 

\une

\subsection{Computations.}

We made some computations to ilustrate how this formulas work.

\enu

\item \emph{Case $n=0$}

Condition (A) implies that $q=e^{-iCx}$.

Condition (B) gives $\|q\|^2$ constant. Hence, $C\in\R$.

\item \emph{Case $n=1$}

Condition (A) gives the Stacionary Non-linear Shrodinger equation:
$$\sigma\|q\|^2q+\frac{1}{2}q_{xx}=\omega q .$$ 

Condition (B) gives an equation in terms of $E$ and $E'$, which may be solved by integration on
an elliptic curve (Riemann Surface of genus $1$):

$$(E')^2=-E^4+4\omega E^2+C_1E+C_2$$

where $C_1$ and $C_2$ are constants.

\item \emph{Case $n=2$}

From equation A, we obtain the equation:

$$iq_{xxx}+6i\sigma\|q\|^2q_x=Cq$$

where $C$ is a constant.

Condition  (B), in this case, is more complicated. 

The equation in terms of E and F is the following:

$$\frac{9}{2}E'E''+\frac{3}{2}EE'''+F'''-\frac{3}{2}E^2F_x-6F_xF-6EFE_x=0$$

\une

\pagebreak

\part{Appendixes}\label{part3}

\begin{appendix}

\section{Proof of Theorem \ref{thm2}}\label{apnA}

We need some definitions and lemmas.

\dfn\label{def1}
Define the following sums
\equ\label{eq3}
S_{l,k}=\sum_{s=k}^{l-k}F_sF_{l-s}'
\uqe
and
\equ\label{eq4}
T_{l,k}=\sum_{s=k}^{l-k}F_sF_{l-s}
\uqe
\nfd

\lmm \label{lem1}
\equ\label{eq5}
\int S_{l,k}=\frac{T_{l,k}}{2}
\uqe
\mml

\begin{proof}
Integrating by parts
$$\int S_{l,k}=\sum_{s=k}^{l-k}\int F_s F_{l-s}'=\sum_{s=k}^{l-k}F_sF_{l-s}-\sum_{s=k}^{l-k}\int F_s'F_{l-s}$$
But,
$$\sum_{s=k}^{l-k}F_s'F_{l-s}=\sum_{s=k}^{l-k} F_s F_{l-s}'=S_{l,k}$$
and, therefore
$$2\int S_{l,k}=\sum_{s=k}^{l-k}F_sF_{l-s}=T_{l,k}.$$
\end{proof}

\lmm \label{lem2}
\equ\label{eq6}
T_{l,k}-T_{l,k+1}=2F_kF_{l-k}
\uqe
\mml

\begin{proof}
Follows from the explicite computation
$$\sum_{s=k}^{l-k}F_sF_{l-s}-\sum_{s=k+1}^{l-k-1}F_sF_{l-s}=F_kF_{l-k}+F_{l-k}F_k.$$
\end{proof}

\lmm\label{lem3}
\equ\label{eq7}
\phi_n^2=\sum_{m=n}^{2n}(4\lambda)^mT_{2n-m-2,-1}+\sum_{m=0}^{n-1}(4\lambda)^mT_{2n-m-2,n-m-1}
\uqe
\mml

\prf
Now, since
$$\phi_n^2=\left(\sum_{r=0}^n (4\lambda)^r F_{n-r-1}\right)\left(\sum_{s=0}^n (4\lambda)^s F_{n-s-1}\right),$$
the terms of degree $m=r+s$ are of the form
$$(4\lambda)^{r+s}F_{n-r-1}F_{n-s-1}$$
and therefore, the sum of the subindexes of $F$ is $2n-m-2$.
Now, since $s=m-r$, we can write this terms only with the variables $m$ and $r$ as
$$(4\lambda)^{m}F_{n-r-1}F_{n-m+r-1},$$
where $n\geq r\geq 0$ or equivalently $n-1\geq n-r-1\geq -1$.
Now, the restriction on the posible values of $r$ given by the degrees in the polynomial $\phi_n$,
is also valid for $s=m-r$. Hence, it must be satisfied also that
$n\geq m-r \geq 0$ and equivalently $n-1\geq n-m+r-1\geq -1$.

The proof of the lemma is divided in two cases.

\enu

\item If $m\geq n$, we list the terms in a descending way from $r=n$ to $r=m-n$.
Notice that when $r=m-n$, the number $n-m+r-1=-1$, its smallest value (or biggest for $s$, $s=n$).
Hence, we have the sum
$$F_{-1}F_{2n-m-1}+F_0F_{2n-m-1}+\cdots+F_{2n-m-1}F_{-1}=T_{2n-m-2,-1}$$

\item If $n>m$, we list the terms in an ascending way from $r=0$ to $m$,
to obtain the sum
$$F_{n-1}F_{n-m-1}+F_{n-2}F_{n-m}+\cdots+F_{n-m-1}F_{-1}=T_{2n-m-2,n-m-1}$$

\une

Both cases prove the formula.

\frp

As in definition (\ref{def1}), but now with the operator $B$, we define the sum

\dfn\label{def2}
\equ\label{eq8}
B_{l,k}=\sum_{s=k}^{l-k}F_sB(F_{l-s})
\uqe
\nfd 

\lmm\label{lem4}
\equ
B_{l,k}=S_{l+1,k}-F_{l+1-k}F_{k}'
\uqe
\mml

\prf
\equ
\begin{split}
B_{l,k} &=\sum_{s=k}^{l-k}F_sB(F_{l-s})=\sum_{s=k}^{l-k}F_sF_{l-s+1}'\\
&=\left(\sum_{s=k}^{l+1-k}F_sF_{l-s+1}'\right)-F_{l+1-k}F_k'=S_{l+1,k}-F_{l+1-k}F_{k}'
\end{split}
\uqe
\frp

Now, we use lema (\ref{lem4}) to prove a formula for $\phi_nB(\phi_n)$ similar to (\ref{eq7}).

\lmm\label{lem5}
\equ\label{eq9}
\phi_nB(\phi_n)=\sum_{m=n}^{2n}(4\lambda)^mS_{2n-m-1,-1}+\sum_{m=0}^{n-1}(4\lambda)^m\left[S_{2n-m-1,n-m-1}-F_{n}F_{n-m-1}'\right]
\uqe
\mml

\prf
First, notice that the same arguments used in the  proof of lemma (\ref{lem3}) can be used to prove
that
\equ\label{eq10}
\phi_nB(\phi_n)=\sum_{m=n}^{2n}(4\lambda)^mB_{2n-m-2,-1}+\sum_{m=0}^{n-1}(4\lambda)^mB_{2n-m-2,n-m-1}
\uqe
But, by lemma (\ref{lem4}), we have
$$B_{2n-m-2,-1}=S_{2n-m-1,-1}-F_{2n-m}F'_{-1}=S_{2n-m-1,-1},$$
(since $F_{-1}$ is constant) and 
$$B_{2n-m-2,n-m-1}=S_{2n-m-1,n-m-1}-F_{n}F_{n-m-1}'.$$
\frp

Now, we are in condition to giving the proof of Theorem (\ref{thm2})

\begin{proof}

Using lemmas (\ref{lem5}) and (\ref{lem1}), we have that
$$\int \phi_nB(\phi_n)dx=\sum_{m=n}^{2n}(4\lambda)^m\frac{T_{2n-m-1,-1}}{2}+\sum_{m=0}^{n-1}(4\lambda)^m\left[\frac{T_{2n-m-1,n-m-1}}{2}-\int F_{n}F_{n-m-1}'dx\right],$$
and multiplying lemma (\ref{lem3}) by $2\lambda$
$$2\lambda\phi_n^2=\sum_{m=n}^{2n}(4\lambda)^{m+1}\frac{T_{2n-m-2,-1}}{2}+\sum_{m=0}^{n-1}(4\lambda)^{m+1}\frac{T_{2n-m-2,n-m-1}}{2}.$$
Substracting
\equ
\begin{split}
\Hc_{n} &=\int\phi_n B(\phi_n) dx -2\lambda\phi_n^2 \\
&=-(4\lambda)^{2n+1}\frac{T_{-2,-1}}{2}+(4\lambda)^{n}\frac{T_{n-1,-1}-T_{n-1,0}}{2}+\\
&\makebox[5mm]{ }+\sum_{m=1}^{n-1}(4\lambda)^m\left[\frac{T_{2n-m-1,n-m-1}-T_{2n-m-1,n-m}}{2}-\int F_{n}F_{n-m-1}'dx\right]\\
&\makebox[5mm]{ }+\left[\frac{T_{2n-1,n-1}}{2}-\int F_nF_{n-1}'\right]
\end{split}
\uqe

Now, using lemma (\ref{lem2}) to compute substractions of T's, and that $T_{-2,-1}=F_{-1}^2$ and $T_{2n-1,n-1}=2F_{n-1}F_n$, we finally conclude

\equ
\begin{split}
\Hc_{n} &=-(4\lambda)^{2n+1}\frac{F_{-1}^2}{2}+(4\lambda)^{n}F_{-1}F_n+\\
&\makebox[5mm]{ }+\sum_{m=1}^{n-1}(4\lambda)^m\left[F_{n-m-1}F_{n}-\int F_{n}F_{n-m-1}'dx\right]\\
&\makebox[5mm]{ }+\left[F_{n-1}F_n-\int F_nF_{n-1}'\right]
\end{split}
\uqe

\end{proof}

\section{One rational symetric identity}\label{apnB}

The main goal of this appendix is to prove the following result. It is important to mention that this theorem appears in \cite{newell} without a proof. The authors here
provided a rigorous proof of this fact. 

\thm\label{mainapn2}
Given a set of different values $\lambda_1$, $\lambda_2$,$\ldots$,$\lambda_{n+1}$, the following identities
are satisfied:
\equ\label{a1}
\sum_{k=1}^{n+1}\frac{\lambda_k^\mu}{\prod_{j\neq k}\left(\lambda_k-\lambda_j\right)}=
\begin{cases}
0 & \mbox{ if }0\leq \mu\leq n-1 \cr
1 & \mbox{ if }\mu=n
\end{cases}
\uqe
\mht

Define the following polynomial of degree $n$ on $x$.
\equ\label{a2}
K_n(x;\lambda_1,\lambda_2,\ldots,\lambda_n)=\prod_{i=1}^n(x-\lambda_i).
\uqe

Now, by induction define the following polynomial in $n+1$ variables.
\equ\label{a3}
P_{n+1}(x_1,x_2,\ldots,x_{n+1})=K_n(x_1;x_2,x_3,\ldots,x_{n+1})P_n(x_2,x_3,\ldots,x_{n+1}),
\uqe
starting with $P_1(x_1)=1$ and $P_0=1$.

These polynomials can also be express by the product formula:
\equ\label{a4}
P_{n+1}(x_1,x_2,\ldots,x_{n+1})=\prod_{1\leq i<j\leq n+1}(x_i-x_j)
\uqe

\exa
Next of these polynomials are:
$$P_2(x_1,x_2)=K_1(x_1;x_2)=x_1-x_2$$
$$P_3(x_1,x_2,x_3)=K_2(x_1;x_2,x_3)P_2(x_2,x_3)=(x_1-x_2)(x_1-x_3)(x_2-x_3)$$
\axe

\lmm\label{a5}
The polynomial $P_{n+1}$ is anti-symetric. (If we interchange any two variables, there is a change
of sign on the polynomial.)
\mml

\begin{proof}
It is sufficient to prove the lemma for two adjacent variables. Now,
given an index $1\leq k\leq n$, we observe that when interchanging the variable $x_k$ with $x_{k+1}$
\enu
\item the terms $x_j-x_k$ and $x_j-x_{k+1}$ are interchanged, for $j<k$.
\item  the terms $x_k-x_j$ and $x_{k+1}-x_j$ are interchanged, for $j>k+1$.
\item the term $x_k-x_{k+1}$ becomes $x_{k+1}-x_k=-(x_k-x_{k+1})$.
\une

Hence, $P_{n+1}$ changes only one sign.

\end{proof}

We need to introduce some notation. If $(x_1,x_2,\ldots,x_n)$ denotes a vector, we will denote
the vector interchange of the variables $x_j$ and $x_k$, with $j<k$ (leaving the other variables
fixed) by $(x_1,\ldots\overleftrightarrow{x_j\ldots x_k}\ldots,x_n)$. If the two variables are adjacent,
we will just denote $(x_1,\ldots\overleftrightarrow{x_j,x_{j+1}}\ldots,x_n)$.

\crl\label{a6}
If $x_i=x_j$ for some pair $i<j$, then $P_n(x_1,x_2,\ldots,x_n)=0$.
\lrc

\begin{proof}
Interchanging  $x_i$ with $x_j$, does not change $P_n(x_1,x_2,\ldots,x_n)$. Hence,
$$P(x_1,\ldots\overleftrightarrow{x_j\ldots x_k}\ldots,x_n)=P_n(x_1,x_2,\ldots,x_n)$$
and also
$$P(x_1,\ldots\overleftrightarrow{x_j\ldots x_k}\ldots,x_n)=-P_n(x_1,x_2,\ldots,x_n)$$ 
\end{proof}

Also, denote $(x_1,\ldots,\widehat{x_i},\ldots,x_n)$, the vector obtained from $(x_1,x_2,\ldots,x_n)$,
after removing the $x_i$ variable.

\dfn
Define the polynomials
\equ\label{a7}
Q_{n+1,j}(x_1,\dots,x_{n+1})=\sum_{i=1}^{n+1}(-1)^{i+1}x_i^jP_n(x_1,\ldots,\widehat{x_i},\ldots,x_{n+1}).
\uqe
\nfd

Observe that, since the polynomial $P_n$ has degree $n-1$ in the first variable $x_1$, the polynomial $Q_{n+1,j}$
has degree at most $n-1$, if $j<n$, and has degree exactly $j$, if $j\geq n$ in the first variable $x_1$.

\lmm\label{a8}
The distinct values $\lambda_2$,$\lambda_3$,$\ldots$,$\lambda_{n+1}$ are roots of the polynomial
$Q_{n+1,j}(X,\lambda_2,\dots,\lambda_{n+1})$.
\mml

\begin{proof}
Evaluating in each $\lambda_k$,
\equ
\begin{split}
Q_{n+1,j}(\lambda_k,\lambda_2,\dots,\lambda_{n+1}) &=\sum_{i=1}^{n+1}(-1)^{i+1}\lambda_i^jP_n(\lambda_k,\ldots,\widehat{\lambda_i},\ldots,\lambda_{n+1})\cr
&= \lambda_k^jP_n(\lambda_2,\lambda_3,\ldots,\lambda_{n+1})+(-1)^{i+1}\lambda_k^jP_n(\lambda_k,\lambda_2,\lambda_3,\ldots,\widehat{\lambda_{k}},\ldots,\lambda_{n+1})
\end{split}
\uqe
since all the terms with two $\lambda_k$ are equal to zero by \ref{a6}.

Now, 
$$P_n(\lambda_k,\lambda_2,\lambda_3,\ldots,\widehat{\lambda_{k}},\ldots,\lambda_{n+1})=(-1)^{k-2} P_n(\lambda_2,\lambda_3,\ldots,\lambda_{n+1}),$$ thus 

\equ
\begin{split}
Q_{n+1,j}(\lambda_k,\lambda_2,\dots,\lambda_{n+1}) &= P_n(\lambda_2,\lambda_3,\ldots,\lambda_{n+1})+(-1)^{2k-1}P_n(\lambda_2,\lambda_3,\ldots,\lambda_{n+1})\cr
&= 0.
\end{split}
\uqe

\end{proof}

\crl\label{a9}
If $j<n$, then
\equ\label{a10}
Q_{n+1,j}(x_1,\dots,x_{n+1})=0.
\uqe
\lrc

\begin{proof}
If $\lambda_2$,$\lambda_3$,$\dots$,$\lambda_{n+1}$ are $n$ distinct values,
$Q_{n+1,j}(x_1,\lambda_2,\lambda_2\dots,\lambda_{n+1})$ is a polynomial in $x_1$ of degree $n-1$
with $n$ roots by \ref{a8}. 

Hence, $Q_{n+1,j}(x_1,\lambda_2,\lambda_2\dots,\lambda_{n+1})=0$.
\end{proof}

A few inductive identities regarding the polynomial $P_n(x_1,\ldots,x_n)$ can be deduced using \ref{a9}.

\lmm\label{a11}
The following identity is true
\equ
P_n(x_1,\ldots,x_n)=\sum_{i=1}^n (-1)^{i+n}x_1\cdots \widehat{x_i}\cdots x_n P_{n-1}(x_1,\ldots,\widehat{x_i},\ldots,x_n)
\uqe
\mml

\begin{proof}
First, notice that
$$P_n(0,x_2,\cdots,x_n)=(-1)^{n-1}x_2\cdots x_n P_{n-1}(x_2,\cdots,x_n).$$
Hence,
\equ
\begin{split}
Q_{n+1,0}(0,x_2,\ldots,x_{n+1}) &= P_n(x_2,\ldots,x_{n+1})+\sum_{i=2}^{n+1}(-1)^{i+1}P_n(0,x_2,\ldots,\widehat{x_i},\ldots, x_{n+1})\cr
&=  P_n(x_2,\ldots,x_{n+1})+\sum_{i=2}^{n+1}(-1)^{i+n}x_2\cdots\widehat{x_i}\cdots x_n P_n(0,x_2\ldots\widehat{x_i}\ldots x_{n+1})\cr
&= 0\cr
\end{split}
\uqe
Therefore, 
$$ P_n(x_2,\ldots,x_{n+1})=\sum_{i=2}^{n+1}(-1)^{i-1+n}x_2\cdots\widehat{x_i}\cdots x_n P_n(0,x_2,\ldots,\widehat{x_i},\ldots, x_{n+1}),$$
which gives the identity of the lemma, shifting by $-1$ the index on the variables.

\end{proof}

\lmm\label{a12}
The following identity is true
\equ
\begin{split}
P_{n}(x_1,\ldots,x_{n})&= Q_{n,n-1}(x_1,\ldots,x_{n})\\
&=\sum_{i=1}^{n} (-1)^{i+1}x_j^{n-1}P_{n-1}(x_1,\ldots,\widehat{x_i},\ldots,x_{n})
\end{split}
\uqe
\mml

\begin{proof}

First, notice that
$$Q_{1,0}(x)=(-1)^2P_0=1=P_1(x)$$
$$Q_{2,1}(x_1,x_2)=(-1)^2x_1P_1(x_2)+(-1)^3x_2P_1(x_1)=x_1-x_2=P_2(x_1,x_2).$$

To symplify the writing of the proof, we will define:

$$\sum_{j=a}^A c_j = 0 \mbox{ if $a>A$ }$$

We will prove the identity for $n+1$, assuming it is valid for $n$.

Using lemma \ref{a11}, we have

\equ\label{a13}
P_{n+1}(x_1,\ldots,x_{n+1}) =\sum_{i=1}^{n+1}x_1\cdots\widehat{x_i}\cdots x_{n+1} P_n(x_1,\ldots, \widehat{x_i},\ldots,x_{n+1})
\uqe

Now, by the induction hypothesis, each term $P_n(x_1,\ldots, \widehat{x_i},\ldots,x_{n+1})$ is equal to

$$
\sum_{j=1}^{i-1}x_j^{n-1}P_{n-1}(x_1,\ldots,\widehat{x_j},\ldots,\widehat{x_i},\ldots,x_n)+
\sum_{j=i+1}^{n+1}x_j^{n-1}P_{n-1}(x_1,\ldots,\widehat{x_i},\ldots,\widehat{x_j},\ldots,x_n)
$$

Combining with (\ref{a13}), we obtain

\equ
P_{n+1}(x_1,\ldots,x_{n+1}) = \sum_{i=1}^{n+1}\sum_{j=1}^{i-1} A_{i,j}+\sum_{i=1}^{n+1}\sum_{j=i+1}^{n+1} B_{i,j}
\uqe

where

$$A_{i,j}=x_j^n (-1)^{n+i+j}x_1\cdots\widehat{x_j}\cdots\widehat{x_i}\cdots x_nP_{n-1}(x_1,\ldots,\widehat{x_j},\ldots,\widehat{x_i},\ldots,x_n)$$

and

$$B_{i,j}=x_j^n (-1)^{n+i+j+1}x_1\cdots\widehat{x_i}\cdots\widehat{x_j}\cdots x_nP_{n-1}(x_1,\ldots,\widehat{x_i},\ldots,\widehat{x_j},\ldots,x_n)$$

Interchanging the order of summation,

\equ\label{a14}
\begin{split}
P_{n+1}(x_1,\ldots,x_{n+1}) &= \sum_{j=1}^{n+1}\sum_{i=j+1}^{n+1} A_{i,j}+\sum_{j=1}^{n+1}\sum_{i=1}^{j-1} B_{i,j}\\
&= \sum_{j=1}^{n+1}\left[\sum_{i=j+1}^{n+1} A_{i,j}+\sum_{i=1}^{j-1} B_{i,j}\right] 
\end{split}
\uqe

Notice that in the first double summation of (\ref{a14}), when $j=n+1$, the term
$$\sum_{i=j+1}^{n+1} A_{i,j}=\sum_{i=n+2}^{n+1} A_{i,j}=0$$
by definition, but this agrees with the fact that in $A_{i,j}$, $n+1=j<i$, 
which gives no term $i$, since there are only $n+1$ variables.
A similar situation occurs in the second summation when $j=1$, then
$$\sum_{i=1}^{j-1} B_{i,j}=\sum_{i=1}^{0} B_{i,j}=0,$$
and now $B_{i,j}$ satisfies that $i<j=1$, but we start labeling variables at $1$.

Now, factoring out $(-1)^{j+1}x_j^n$ from the terms  $A_{i,j}$ and  $B_{i,j}$,
we obtain

\equ\label{a15}
P_{n+1}(x_1,\ldots,x_{n+1}) = \sum_{j=1}^{n+1} (-1)^{j+1}x_j^n C_{i,j}
\uqe

with

\equ
\begin{split}
C_{i,j} &=\sum_{i=j+1}^{n+1}(-1)^{n+i-1}x_1\cdots\widehat{x_j}\cdots\widehat{x_i}\cdots x_nP_{n-1}(x_1,\ldots\widehat{x_j}\ldots\widehat{x_i}\ldots,x_n)+\\
&+\sum_{i=1}^{j-1}(-1)^{n+i}x_1\cdots\widehat{x_i}\cdots\widehat{x_j}\cdots x_nP_{n-1}(x_1,\ldots\widehat{x_i}\ldots\widehat{x_j}\ldots,x_n).
\end{split}
\uqe

But, lemma \ref{a11} gives

$$C_{i,j}=P_n(x_1,\ldots,\widehat{x_j},\ldots,x_{n+1}).$$

Hence,

\equ
\begin{split}
P_{n+1}(x_1,\ldots,x_{n+1}) &=\sum_{j=1}^{n+1} (-1)^{j+1}x_j^nP_n(x_1,\ldots,\widehat{x_j},\ldots,x_{n+1})\\
&= Q_{n+1,n}(x_1,\ldots,x_{n+1}).
\end{split}
\uqe

\end{proof}

Now, we are ready to proof Theorem \ref{mainapn2}.

\begin{proof}
First, notice
$$K_n(\overleftrightarrow{\lambda_1\ldots \lambda_k},\ldots,\lambda_{n+1})=\prod_{j\neq k}\left(\lambda_k-\lambda_j\right).$$
Hence,

\equ
\label{a15}
\sum_{k=1}^{n+1}\frac{\lambda_k^\mu}{\prod_{j\neq k}\left(\lambda_k-\lambda_j\right)}=\sum_{k=1}^{n+1}\frac{\lambda_k^\mu}{K_n(\overleftrightarrow{\lambda_1\ldots \lambda_k},\ldots,\lambda_{n+1})}.
\uqe

Now, permuting $x_1$ with $x_k$ in equation (\ref{a3}) we obtain

$$
P_{n+1}(\overleftrightarrow{x_1\ldots x_k}\ldots,x_{n+1})=K_n(\overleftrightarrow{x_1\ldots x_k}\ldots,x_{n+1})P_n(x_2\ldots,x_{k-1},x_1,x_{k+1}\ldots,x_{n+1})
$$

But, now moving the variable $x_1$ at the begining in the $P_n$ term on the right of the above equation, we have
$$P_n(x_2\ldots,x_{k-1},x_1,x_{k+1}\ldots,x_{n+1})=(-1)^{k-2}P_n(x_1,\ldots\widehat{x_k}\ldots,x_n)$$

And, using that
$$P_{n+1}(\overleftrightarrow{x_1\ldots x_k}\ldots,x_{n+1})=-P_{n+1}(x_1,\ldots,x_{n+1}),$$
we obtain
$$
P_{n+1}(x_1,\ldots,x_{n+1})=(-1)^{k-1}K_n(\overleftrightarrow{x_1\ldots x_k}\ldots,x_{n+1})P_n(x_1,\ldots\widehat{x_k}\ldots,x_n).
$$

Hence, equation (\ref{a15}) becomes

\equ
\label{a16}
\begin{split}
\sum_{k=1}^{n+1}\frac{\lambda_k^\mu}{\prod_{j\neq k}\left(\lambda_k-\lambda_j\right)} &=\frac{\sum_{k=1}^{n+1}(-1)^{k-1}\lambda_k^\mu P_n(\lambda_1,\ldots\widehat{\lambda_k}\ldots,\lambda_n)}{P_{n+1}(\lambda_1,\ldots,\lambda_{n+1})}\\
&=\frac{Q_{n+1,\mu}(\lambda_1,\ldots,\lambda_{n+1})}{P_{n+1}(\lambda_1,\ldots,\lambda_{n+1})}
\end{split}
\uqe

Therefore, the above sum is $0$ if $0\leq \mu \leq n-1$ (corollary \ref{a9})
and $1$ if $\mu=n$ (lemma \ref{a12}).

\end{proof}

\section{Note on succesive derivatives of $E=\frac{iq_x}{q}$}\label{apnC}

The first derivative of $E$ is
$$E'=\frac{iq_{xx}}{q}-i\frac{q_x^2}{q^2}$$
Setting 
\equ\label{E2:definition}
E_{(2)}=\frac{iq_{xx}}{q},
\uqe
 we can write the first derivative as:
\equ\label{E:first:derivative}
E'=E_{(2)}+iE^2
\uqe

Now, define:
\equ\label{En:definition}
E_{(n)}=\frac{iq^{(n)}}{q}.
\uqe
Hence, we easily compute:
\equ\label{En:derivative}
E_{(n)}'=E_{(n+1)}+iE_{(n)}E
\uqe

Using this notation, we can easily compute derivatives of $E$ of superior order.
For example,
\equ\label{E:second:derivative}
\begin{split}
E'' &=E_{(2)}'+2iE\cdot E'\\
&=E_{(3)}+iE_{(2)}E+2iE\cdot\left(E_{(2)}+iE^2\right)\\
&=E_{(3)}+3iE_{(2)}E-2E^3
\end{split}
\uqe

and

\equ
\begin{split}\label{E:third:derivative}
E''' &=E_{(3)}'+3i\left(E_{(2)}E\right)'-6E^2E'\\
&=E_{(4)}+iE_{(3)}E+3i\left[E_{(2)}'E+E_{(2)}E'\right]-6E^2\cdot\left[E_{(2)}+iE^2\right]\\
&=E_{(4)}+iE_{(3)}E-6E^2E_{(2)}-6iE^4+3iE_{(2)}'E+3iE_{(2)}E'\\
&=E_{(4)}+iE_{(3)}E-6E^2E_{(2)}-6iE^4+3iE_{(2)}'E+3iE_{(2)}\left(E_{(2)}+iE^2\right)\\
&=E_{(4)}+iE_{(3)}E-9E^2E_{(2)}-6iE^4+3iE_{(2)}^2+3iE_{(2)}'E\\
&=E_{(4)}+iE_{(3)}E-9E^2E_{(2)}-6iE^4+3iE_{(2)}^2+3iE\left(E_{(3)}+iE_{(2)}E\right)\\
&=E_{(4)}+4iE_{(3)}E-12E^2E_{(2)}-6iE^4+3iE_{(2)}^2
\end{split}
\uqe

\end{appendix}






\begin{thebibliography}{99}



\bibitem{scott} A.S. Scott. {\em The development of nonlinear Science.}  Rivista del Nuovo Cimento, vol. \textbf{27}, Issue 10, p.1-115.

\bibitem{KdV_seminal}  Korteweg, D. J.; de Vries, G. (1895), {\em On the Change of Form of Long Waves Advancing in a Rectangular Canal, and on a New Type of Long Stationary Waves}, Philosophical Magazine, \textbf{39} (240): 422?443,

\bibitem{dubrovin_1} B.A. Dubrovin. {\em Theta functions and non-linear equations.} Russian Math. Surveys \textbf{36}:2 (1981), 11-92.

\bibitem{krichever_1} I.M. Krichever. {\em Methods of algebraic geometry in the theory of non-linear equations.} Russian Math. Surveys \textbf{32}:6 (1977), 185-213. 

\bibitem{krichever_novikov}   I.M. Krichever and S.P. Novikov.  {\em   Integration of non-linear equations by methods of algebraic geometry.} Functional Anal. Appl. \textbf{11}:1 (1977), 12-26.

\bibitem{dubrovin_matveev_novikov}  B.A. Dubrovin, V.B. Matveev and S.P. Novikov. {\em Non-linear equations of Korteweg-deVries type, finite-zone operators, and Abelian varieties}. Russian Math. Surveys \textbf{31}:1 (1976), 59-146

\bibitem{} S.P. Novikov, {\em A periodic problem for the Korteweg-de Vries equation, I.} Functional Anal. Appl. \textbf{8} (1974), 236-246.

\bibitem{} B.A. Dubrovin and S.P. Novikov. {\em A periodic problem for the Korteweg-de Vries and Sturm-Iiouville equations. Their connection with algebraic geometry.} Soviet Math. Dokl. \textbf{15} (1974), 1597-1601.

\bibitem{} I.M. Krichever. {\em Algebraic curves and non-linear difference equations.}  Russian Math. Survey \textbf{33}:4 (1978), 255-256.

\bibitem{} B.A. Dubrovin and S.P. Novikov.{\em Periodic and conditionally periodic analogues of the multi-soliton solutions of the Korteweg-de Vries equation.} Soviet Physics JETP \textbf{40} (1974), 1058-1063.

\bibitem{}  B.A. Dubrovin . {\em Periodic problems for the Korteweg-de Vries equation in the class of finite-zone potentials.} Functional Anal. Appl. \textbf{9} (1975), 215-223.

\bibitem{} I.M. Krichever. {\em An algebraic-geometric construction of the Zakharov-Shabat equations and their periodic solutions.}  Soviet Math. Dokl. \textbf{17} (1976), 394-397.

\bibitem{} A.P. Its and V.B. Matveev. {\em Schr\"odinger operators with a finite-zone spectrum and the jV-soliton solutions of the Korteweg-de Vries equation.} Teoret. Mat. Fiz. \textbf{23}
(1975), 51-67.

\bibitem{flaschka_1}  H. Flaschka. {\em Personal communication.}

\bibitem{newell} A.C. Newell. {\em Solitons in Mathematics and Physics.} Society for the Industrial and Applied Mathematics. Philadelphia, PA. USA. 1985.

\bibitem{zajarov} V.E. Zakharov,  and A.B. Shabat. {\em Exact Theory of Two Dimensional Self-Focusing and One-Dimensional Self-Modulation of Waves in Nonlinear Media.} Soviet Physics-JETP. Vol. \textbf{ 34}, No. 1, (January 1972) 62-69.

\bibitem{Kamchatnov_Kraenkel_Umarov_1}  A.M. Kamchatnov, and R.A. Kraenkel, and B.A. Umarov, Phys. Lett. A 287, 223  2001.

\bibitem{Kamchatnov_Kraenkel} A.M. Kamchatnov and R.A. Kraenkel, J. Phys. A 35, L13  2002 .

\bibitem{Kamchatnov_Kraenkel_Umarov_2} A.M. Kamchatnov, R. A. Kraenkel, and B. A. Umarov. {\em Asymptotic soliton train solutions of the defocusing nonlinear Schr\"odinger equation}. Physical Review E. \textbf{66}, 036609,  (2002).

\bibitem{Kamchatnov_1} A.M. Kamchatnov. {\em New approach to periodic solutions of integrable equations and nonlinear theory of modulational instability}. Physics Reports \textbf{286} (1997) 199-270

\bibitem{Kamchatnov_2} A.M. Kamchatnov. {\em Nonlinear periodic waves and their modulations: an introductory course.} World Scientific Publishing. Singapore. 2000.

\bibitem{espinola_portillo} Esp\'{\i}nola-Rocha, J.A. \& Portillo-Bobadilla, F.X. (2017) {\it Multiplicative Operators in the spectral problem of integrable systems.}  Journal of Integrable Systems (To be submited).

  
\bibitem{kdv} Korteweg, D. J.; de Vries, G. (1895), "On the Change of Form of Long Waves Advancing in a Rectangular Canal, and on a New Type of Long Stationary Waves", Philosophical Magazine, 39 (240): 422?443,

\bibitem{ggkm} C. Gardner, J. Greene, M. Kruskal and R. Miura, {\em A Method for Solving the Korteweg-DeVries Equation.} Phys. Rev. Letters, \textbf{ 19} (1967) 1095.

\bibitem{drazin}  P.G. Drazin and R.S. Johnson, {\em Solitons: an introduction.} Cambridge University Press, 1989.

\bibitem{novikov} S. Novikov,   S.V. Manakov,  L.P. Pitaevskii, and V.E. Zakharov, {\em Theory of Solitons: The Inverse Scattering Method.} New York : Consultants Bureau, 1984.

\bibitem{zakharov_faddeev} V.E. Zakharov and L.D. Faddeev. {\em Korteweg-de Vries equation: A completely integrable Hamiltonian system.} Functional analysis and its applications \textbf{5} (4), 280-287.

\bibitem{lax_1}  P. Lax, \textit{Integrals of nonlinear equations of evolutions and solitary waves}. Communications on Pure and Applied Mathematics. \textbf{XXI}. 467-490 (1968).

\bibitem{lax_2}   P. Lax, \textit{Periodic solutions of the KdV equation}. Communications on Pure and Applied Mathematics. \textbf{XXVIII}. 141-188 (1975).


\end{thebibliography}
\end{document}